\catcode`\@=11

\font\tenmsa=msam10
\font\sevenmsa=msam7
\font\fivemsa=msam5
\font\tenmsb=msbm10
\font\sevenmsb=msbm7
\font\fivemsb=msbm5
\newfam\msafam
\newfam\msbfam
\textfont\msafam=\tenmsa  \scriptfont\msafam=\sevenmsa
  \scriptscriptfont\msafam=\fivemsa
\textfont\msbfam=\tenmsb  \scriptfont\msbfam=\sevenmsb
  \scriptscriptfont\msbfam=\fivemsb

\def\hexnumber@#1{\ifnum#1<10 \number#1\else
 \ifnum#1=10 A\else\ifnum#1=11 B\else\ifnum#1=12 C\else
 \ifnum#1=13 D\else\ifnum#1=14 E\else\ifnum#1=15 F\fi\fi\fi\fi\fi\fi\fi}

\def\msa@{\hexnumber@\msafam}
\def\msb@{\hexnumber@\msbfam}
\mathchardef\boxdot="2\msa@00
\mathchardef\boxplus="2\msa@01
\mathchardef\boxtimes="2\msa@02
\mathchardef\square="0\msa@03
\mathchardef\blacksquare="0\msa@04
\mathchardef\centerdot="2\msa@05
\mathchardef\lozenge="0\msa@06
\mathchardef\blacklozenge="0\msa@07
\mathchardef\circlearrowright="3\msa@08
\mathchardef\circlearrowleft="3\msa@09
\mathchardef\rightleftharpoons="3\msa@0A
\mathchardef\leftrightharpoons="3\msa@0B
\mathchardef\boxminus="2\msa@0C
\mathchardef\Vdash="3\msa@0D
\mathchardef\Vvdash="3\msa@0E
\mathchardef\vDash="3\msa@0F
\mathchardef\twoheadrightarrow="3\msa@10
\mathchardef\twoheadleftarrow="3\msa@11
\mathchardef\leftleftarrows="3\msa@12
\mathchardef\rightrightarrows="3\msa@13
\mathchardef\upuparrows="3\msa@14
\mathchardef\downdownarrows="3\msa@15
\mathchardef\upharpoonright="3\msa@16

\mathchardef\downharpoonright="3\msa@17
\mathchardef\upharpoonleft="3\msa@18
\mathchardef\downharpoonleft="3\msa@19
\mathchardef\rightarrowtail="3\msa@1A
\mathchardef\leftarrowtail="3\msa@1B
\mathchardef\leftrightarrows="3\msa@1C
\mathchardef\rightleftarrows="3\msa@1D
\mathchardef\Lsh="3\msa@1E
\mathchardef\Rsh="3\msa@1F
\mathchardef\rightsquigarrow="3\msa@20
\mathchardef\leftrightsquigarrow="3\msa@21
\mathchardef\looparrowleft="3\msa@22
\mathchardef\looparrowright="3\msa@23
\mathchardef\circeq="3\msa@24
\mathchardef\succsim="3\msa@25
\mathchardef\gtrsim="3\msa@26
\mathchardef\gtrapprox="3\msa@27
\mathchardef\multimap="3\msa@28
\mathchardef\therefore="3\msa@29
\mathchardef\because="3\msa@2A
\mathchardef\doteqdot="3\msa@2B

\mathchardef\triangleq="3\msa@2C
\mathchardef\precsim="3\msa@2D
\mathchardef\lesssim="3\msa@2E
\mathchardef\lessapprox="3\msa@2F
\mathchardef\eqslantless="3\msa@30
\mathchardef\eqslantgtr="3\msa@31
\mathchardef\curlyeqprec="3\msa@32
\mathchardef\curlyeqsucc="3\msa@33
\mathchardef\preccurlyeq="3\msa@34
\mathchardef\leqq="3\msa@35
\mathchardef\leqslant="3\msa@36
\mathchardef\lessgtr="3\msa@37
\mathchardef\backprime="0\msa@38
\mathchardef\risingdotseq="3\msa@3A
\mathchardef\fallingdotseq="3\msa@3B
\mathchardef\succcurlyeq="3\msa@3C
\mathchardef\geqq="3\msa@3D
\mathchardef\geqslant="3\msa@3E
\mathchardef\gtrless="3\msa@3F
\mathchardef\sqsubset="3\msa@40
\mathchardef\sqsupset="3\msa@41
\mathchardef\trianglerighteq="3\msa@44
\mathchardef\trianglelefteq="3\msa@45
\mathchardef\bigstar="0\msa@46
\mathchardef\between="3\msa@47
\mathchardef\blacktriangledown="0\msa@48
\mathchardef\blacktriangleright="3\msa@49
\mathchardef\blacktriangleleft="3\msa@4A
\mathchardef\blacktriangle="0\msa@4E
\mathchardef\triangledown="0\msa@4F
\mathchardef\eqcirc="3\msa@50
\mathchardef\lesseqgtr="3\msa@51
\mathchardef\gtreqless="3\msa@52
\mathchardef\lesseqqgtr="3\msa@53
\mathchardef\gtreqqless="3\msa@54
\mathchardef\Rrightarrow="3\msa@56
\mathchardef\Lleftarrow="3\msa@57
\mathchardef\veebar="2\msa@59
\mathchardef\barwedge="2\msa@5A
\mathchardef\doublebarwedge="2\msa@5B
\mathchardef\angle="0\msa@5C
\mathchardef\measuredangle="0\msa@5D
\mathchardef\sphericalangle="0\msa@5E
\mathchardef\varpropto="3\msa@5F
\mathchardef\smallsmile="3\msa@60
\mathchardef\smallfrown="3\msa@61
\mathchardef\Subset="3\msa@62
\mathchardef\Supset="3\msa@63
\mathchardef\Cup="2\msa@64

\mathchardef\Cap="2\msa@65

\mathchardef\curlywedge="2\msa@66
\mathchardef\curlyvee="2\msa@67
\mathchardef\leftthreetimes="2\msa@68
\mathchardef\rightthreetimes="2\msa@69
\mathchardef\subseteqq="3\msa@6A
\mathchardef\supseteqq="3\msa@6B
\mathchardef\bumpeq="3\msa@6C
\mathchardef\Bumpeq="3\msa@6D
\mathchardef\lll="3\msa@6E

\mathchardef\ggg="3\msa@6F

\mathchardef\circledS="0\msa@73
\mathchardef\pitchfork="3\msa@74
\mathchardef\dotplus="2\msa@75
\mathchardef\backsim="3\msa@76
\mathchardef\backsimeq="3\msa@77
\mathchardef\complement="0\msa@7B
\mathchardef\intercal="2\msa@7C
\mathchardef\circledcirc="2\msa@7D
\mathchardef\circledast="2\msa@7E
\mathchardef\circleddash="2\msa@7F
\def\ulcorner{\delimiter"4\msa@70\msa@70 }
\def\urcorner{\delimiter"5\msa@71\msa@71 }
\def\llcorner{\delimiter"4\msa@78\msa@78 }
\def\lrcorner{\delimiter"5\msa@79\msa@79 }
\def\yen{\mathhexbox\msa@55 }
\def\checkmark{\mathhexbox\msa@58 }
\def\circledR{\mathhexbox\msa@72 }
\def\maltese{\mathhexbox\msa@7A }
\mathchardef\lvertneqq="3\msb@00
\mathchardef\gvertneqq="3\msb@01
\mathchardef\nleq="3\msb@02
\mathchardef\ngeq="3\msb@03
\mathchardef\nless="3\msb@04
\mathchardef\ngtr="3\msb@05
\mathchardef\nprec="3\msb@06
\mathchardef\nsucc="3\msb@07
\mathchardef\lneqq="3\msb@08
\mathchardef\gneqq="3\msb@09
\mathchardef\nleqslant="3\msb@0A
\mathchardef\ngeqslant="3\msb@0B
\mathchardef\lneq="3\msb@0C
\mathchardef\gneq="3\msb@0D
\mathchardef\npreceq="3\msb@0E
\mathchardef\nsucceq="3\msb@0F
\mathchardef\precnsim="3\msb@10
\mathchardef\succnsim="3\msb@11
\mathchardef\lnsim="3\msb@12
\mathchardef\gnsim="3\msb@13
\mathchardef\nleqq="3\msb@14
\mathchardef\ngeqq="3\msb@15
\mathchardef\precneqq="3\msb@16
\mathchardef\succneqq="3\msb@17
\mathchardef\precnapprox="3\msb@18
\mathchardef\succnapprox="3\msb@19
\mathchardef\lnapprox="3\msb@1A
\mathchardef\gnapprox="3\msb@1B
\mathchardef\nsim="3\msb@1C
\mathchardef\napprox="3\msb@1D
\mathchardef\nsubseteqq="3\msb@22
\mathchardef\nsupseteqq="3\msb@23
\mathchardef\subsetneqq="3\msb@24
\mathchardef\supsetneqq="3\msb@25
\mathchardef\subsetneq="3\msb@28
\mathchardef\supsetneq="3\msb@29
\mathchardef\nsubseteq="3\msb@2A
\mathchardef\nsupseteq="3\msb@2B
\mathchardef\nparallel="3\msb@2C
\mathchardef\nmid="3\msb@2D
\mathchardef\nshortmid="3\msb@2E
\mathchardef\nshortparallel="3\msb@2F
\mathchardef\nvdash="3\msb@30
\mathchardef\nVdash="3\msb@31
\mathchardef\nvDash="3\msb@32
\mathchardef\nVDash="3\msb@33
\mathchardef\ntrianglerighteq="3\msb@34
\mathchardef\ntrianglelefteq="3\msb@35
\mathchardef\ntriangleleft="3\msb@36
\mathchardef\ntriangleright="3\msb@37
\mathchardef\nleftarrow="3\msb@38
\mathchardef\nrightarrow="3\msb@39
\mathchardef\nLeftarrow="3\msb@3A
\mathchardef\nRightarrow="3\msb@3B
\mathchardef\nLeftrightarrow="3\msb@3C
\mathchardef\nleftrightarrow="3\msb@3D
\mathchardef\divideontimes="2\msb@3E
\mathchardef\varnothing="0\msb@3F
\mathchardef\nexists="0\msb@40
\mathchardef\mho="0\msb@66
\mathchardef\thorn="0\msb@67
\mathchardef\beth="0\msb@69
\mathchardef\gimel="0\msb@6A
\mathchardef\daleth="0\msb@6B
\mathchardef\lessdot="3\msb@6C
\mathchardef\gtrdot="3\msb@6D
\mathchardef\ltimes="2\msb@6E
\mathchardef\rtimes="2\msb@6F
\mathchardef\shortmid="3\msb@70
\mathchardef\shortparallel="3\msb@71
\mathchardef\smallsetminus="2\msb@72
\mathchardef\thicksim="3\msb@73
\mathchardef\thickapprox="3\msb@74
\mathchardef\approxeq="3\msb@75
\mathchardef\succapprox="3\msb@76
\mathchardef\precapprox="3\msb@77
\mathchardef\curvearrowleft="3\msb@78
\mathchardef\curvearrowright="3\msb@79
\mathchardef\digamma="0\msb@7A
\mathchardef\varkappa="0\msb@7B
\mathchardef\hslash="0\msb@7D
\mathchardef\hbar="0\msb@7E
\mathchardef\backepsilon="3\msb@7F
\def\Bbb{\ifmmode\let\next\Bbb@\else
 \def\next{\errmessage{Use \string\Bbb\space only in math mode}}\fi\next}
\def\Bbb@#1{{\Bbb@@{#1}}}
\def\Bbb@@#1{\fam\msbfam#1}

\catcode`\@=\active

\def\CL{\hbox{{$\cal L$}}}

\def\CR{\hbox{{$\cal R$}}} 
 
\def\CM{\hbox{{$\cal M$}}}

\def\CQ{\hbox{{$\cal Q$}}}  \def\CE{\hbox{{$\cal E$}}}

\def\cg{{\sl g}}
\def\R{{\Bbb R}}
\def\C{{\Bbb C}}
\def\Z{{\Bbb Z}}

\def\cu{{\sl u}}

\def\lform{\hbox{$\sqcup$}\llap{\hbox{$\sqcap$}}}
\def\h{{{1\over2}}}
\def\eps{{\epsilon}}

\def\swap{{\leftrightarrow}}

\def\<{\langle}
\def\>{\rangle}

\def\lcross{{>\!\!\!\triangleleft}}

\def\dcross{{\bowtie}}
\def\codcross{{\blacktriangleright\!\!\blacktriangleleft}}
\def\lbiprod{{\lcross}} %%% temp

\def\tens{\mathop{\otimes}}
\def\la{{\triangleright}}
\def\isom{{\cong}}

\def\Ad{{\rm Ad}}

\def\ev{{\rm ev}}
\def\coev{{\rm coev}}

\def\id{{\rm id}}

\def\dila{{\varsigma}}

\def\o{{}_{(1)}}\def\t{{}_{(2)}}\def\th{{}_{(3)}}\def\fo{{}_{(4)}}
\def\bo{{}^{\bar{(1)}}}\def\bt{{}^{\bar{(2)}}}
\def\tildebo{{}^{\tilde{(1)}}}\def\tildebt{{}^{\tilde{(2)}}}
\def\uo{{{}^{(1)}}}\def\ut{{{}^{(2)}}}\def\uth{{{}^{(3)}}}
\def\Bo{{{}_{\underline{(1)}}}}\def\Bt{{{}_{\underline{(2)}}}}
\def\Bth{{{}_{\underline{(3)}}}}

\def\vb{{\sl e}}\def\cvb{{\sl f}}\def\qg{{H}}\def\vecu{{\bf u}}

\def\extd{{\rm d}}

\def\proof{\goodbreak\noindent{\bf Proof\quad}}
\def\endproof{{\ $\lform$}\bigskip }

\def\text#1{{\rm #1}}
\def\note#1{}
\def\und#1{{\underline #1}}

\def\nquad{{\!\!\!\!\!\!}}

\def\equad{\nquad}
\def\eqn#1#2{\begin{equation}#2\label{#1}\end{equation}}
\def\align#1{\begin{eqnarray*}#1\end{eqnarray*}}
\documentstyle[11pt,epsf]{article}
\newtheorem{lemma}{Lemma}[section]
\newtheorem{propos}[lemma]{Proposition}

\newtheorem{corol}[lemma]{Corollary}
\newtheorem{defin}[lemma]{Definition}

\begin{document}
\baselineskip 23pt

{\ }\hskip 4.7in   Damtp/97-73
\vspace{.2in}

\begin{center} {\Large QUANTUM AND BRAIDED GROUP RIEMANNIAN GEOMETRY}
\\ \baselineskip 13pt{\ }
{\ }\\ Shahn Majid\footnote{Royal Society University Research
Fellow and Fellow of Pembroke College, Cambridge, England.}\\ {\
}\\ Department of Applied Mathematics and Theoretical Physics\\
University of Cambridge, Cambridge CB3 9EW, UK\\
www.damtp.cam.ac.uk/user/majid
\end{center}

\begin{center}
August, 1997
\end{center}
\vspace{10pt}
\begin{quote}\baselineskip 13pt
\noindent{\bf ABSTRACT}
We formulate quantum group Riemannian geometry as a gauge theory of
quantum differential forms. We first develop (and slightly
generalise) classical Riemannian geometry in a self-dual manner as
a principal bundle frame resolution and a dual pair of canonical
forms. The role of Levi-Civita connection is naturally generalised
to connections with vanishing torsion and cotorsion, which we
introduce. We then provide the corresponding quantum group and
braided group formulations with the universal quantum differential
calculus. We also give general constructions for examples,
including quantum spheres and quantum planes.
\end{quote}
\baselineskip 23pt

\section{INTRODUCTION}

One of the long-term motivations for noncommutative geometry is a
theory of Riemannian geometry and gravity powerful enough to not
break down in the quantum domain. The quantum groups approach to
such `Planck scale geometry' is initiated in \cite{Ma:pla} in the
context of quantum phases spaces which are quantum groups. The
axioms of a quantum group provide a way to unify the
noncommutativity due to quantisation and the noncocommutativity or
nonAbelianness of the coproduct due to curvature.  To extend  this
programme  to more realistic models, however, one needs a more
general framework not limited to analogues of group manifolds. One
may expect that quantum groups will still play a role as symmetries
or quantum gauge groups and that they will ensure a plentiful
supply of examples such as q-deformed homogeneous spaces. A further
and more general motivation is the q-deformation of geometric
structures to provide a new  regularisation parameter `q' whether
or not it is related to Planck's constant. In some settings, $q$ in
fact generalises the role of $-1$ in bose-fermi statistics, so that
this quantum geometry effectively generalises supergeometry as well
\cite{Ma:introp}. We refer to the text \cite{Ma:book} for much of
this background.

Motivated in part by such considerations, quantum group gauge
theory has already been introduced in \cite{BrzMa:gau} and several
follow-up papers
\cite{Haj:str}\cite{Ma:war95}\cite{Brz:tra}\cite{BrzMa:dif}. It has
also been extended to more general settings
\cite{BrzMa:coa}\cite{Ma:diag}. On the other hand, Riemannian
geometry has not yet been included in this set-up. In the present
paper we provide a start to this programme by giving a fresh
formulation of classical Riemannian geometry as a gauge theory (in
Section~2) and proceeding (in Section~3) to the parallel
definitions in the noncommutative or `quantum' case.  We limit our
results to the `free' or universal differential calculus
$\Omega^1M$ associated to any (not necessarily commutative algebra
$M$), but the case of nonuniversal differential calculus can be
treated similarly and will be developed in detail elsewhere.

Riemannian geometry as a gauge theory is obviously very well known,
in the physics literature under the heading of `vierbeins' and in
the mathematics literature as the theory of so-called `linear' or
frame connections\cite{KobNom:fou}. Our main innovation in the
classical theory is an emphasis on differential forms in place of
vector fields and the use of a general Lie group in place of the
group of linear transformations of local bases. We recall that
usually one views the tangent bundle of a manifold as associated to
the bundle of linear frames (a certain $GL_n$ principal bundle),
and the Levi-Civita connection as induced from a gauge connection
on it. Key to this correspondence is the canonical one-form
$\theta$ on the frame bundle. Recently, in the appendix of
\cite{Haj:str}, it is observed that this form $\theta$ fully
characterises the frame bundle, a result which is presented as `a
theorem which should have been proven 30 years ago'. We extend this
observation now by showing that other constructions in Riemannian
geometry also factor through this canonical 1-form $\theta$. As
with the appendix of \cite{Haj:str}, such results are surely
`known' in some form, in the sense that they implicitly underly the
standard treatments such as\cite{KobNom:fou}, but  they are usually
done for the particular standard $GL_n$  frame bundle. One also
finds separate calculations for $O_n$ and affine frame bundles. By
contrast we provide derivations of the torsion and curvature
tensors (the main formulae of Riemannian geometry) associated to a
any {\em frame resolution} $(P,G,V,\theta)$ (where $P$ is a
principal $G$-bundle and $V$ a representation of $G$), independent
of its detailed form. We have not found such a non-specific
treatment elsewhere, although \cite{KMS:nat} comes close in some
points (notably our Proposition~2.3). As an application,
Proposition~2.5 expresses the torsion tensor as the difference
$T=\extd-\nabla$ between the covariantised exterior derivative on
1-forms and the usual exterior derivative. Moreover, we show that
dual frame resolutions $(P,G,V^*,\gamma)$ is equivalent to a
(not-necessarily symmetric) metric tensor, which is a new
`self-dual' way of thinking about Riemannian geometry. The dual
theory has the roles of $V,\theta$ and $V^*,\gamma$ interchanged,
which also transposes the metric (so the usual symmetric metric is
the self-dual case). We show that the vanishing of the {\em
cotorsion} (the torsion in the dual frame resolution) is a skew
version of metric compatibility. This gives a natural
generalisation of Levi-Civita connections as ones with zero torsion
and zero cotorsion.

Proceeding to the quantum group case, the main problem of quantum
group Riemannian geometry is what quantum group? Surely not
$SO_q(n)$ or $GL_q(n)$ in general as `quantum spaces' are not in
general `manifolds' based on patching $\R_q^n$. In general they are
algebras without necessarily even a dimension $n$. Based on our
formulation of the  classical theory in Section~2, we propose to
consider as `manifold' a possibly noncommutative algebra $M$
equipped with a differential calculus $\Omega^1(M)$ (or higher
forms) as in quantum group gauge theory, and to consider any
`quantum frame resolution' of $(M,\Omega^1(M))$. This is any
quantum principal bundle over $M$ equipped with certain further
structure. In this way, we are  not forced to choose which quantum
group will play the role of structure group of the frame bundle.
Any choice of resolution will serve the purpose of expressing the
quantum (co)tangent bundle as an associated vector bundle and allow
the first steps of quantum group Riemannian geometry as an
application of the gauge theory of associated bundles already in
\cite{BrzMa:gau}. The second and related problem is that one cannot
expect the quantum metric to be symmetric (it typically has some
form or q-symmetry). Again motivated by our formulation of the
classical theory, we drop such a consideration and formulate the
metric a dual frame resolution.

As a modest first application of the quantum group version, we show
in Section~4 that the well-known `parallelisation isomorphism'
$H\tens
\ker\eps\isom
\Omega^1H$  used in the theory of bicovariant differential
calculi on a quantum group $H$ can be understood `geometrically' as
induced by a frame bundle resolution. We also show
(Proposition~4.2) how Fourier transformation on a
finite-dimensional Hopf algebra can be interpreted as inducing a
quantum metric on it. We then turn to more general classes of
examples, including a frame resolution for any quantum homogeneous
space (such a the $q$-sphere) obtained as a homogeneous bundle
$P\to H$, and a frame bundle resolution for any braided group $B$
(such as the quantum plane $M=\R_q^n$) obtained as its
bosonisation. The structure group in the quantum plane example is
$GL_q(n)$, i.e. this is the natural `flat space' example.

The paper concludes with an appendix detailing the corresponding
theory with braid statistics, proven diagrammatically. This is for
completeness as a supplement the diagrammatic or braided group
gauge theory in \cite{BrzMa:coa}\cite{Ma:diag}. The still more
general coalgebra bundle case\cite{BrzMa:coa} is a direction for
further work.

Finally, although we will indicate  briefly how the results
generalise to non-universal differential calculi, the detailed
theory in this case must await the theory of associated bundles
with nonuniversal calculus. The first step, namely the construction
of natural nonuniversal calculi $\Omega^1(P)$ from $\Omega^1(M)$
and nonuniversal calculi on the fiber, has recently been obtained
in \cite{BrzMa:dif}. However, the extension of this to associated
bundles is a separate and quite involved project, to be considered
in a sequel.

Our quantum groups approach should ultimately link up with other
approaches to noncommutative geometry, notably with that of A.
Connes\cite{Con:non}. There one considers only vector bundles (as
projective modules) and abstractly-defined covariant derivations,
and not principal bundles and connection forms (which really need
the quantum groups approach\cite{BrzMa:gau}). Some initial results
relating at least gauge theory  in the two approaches are in
\cite{HajMa:pro}. Among other papers,  \cite{HecSch:lev} have
recently given examples in the vector bundle and `covariant
derivation' approach with the base $M$ a standard quantum
group such as $SO_q(n)$ (with nonuniversal differential calculus);
it should ultimately be possible to understand this as part of a
general theory of the type presented here. Finally, another
attempt at frame bundles can be found in \cite{Dur:fra}. By
contrast, we have adopted the proposal for frame bundles in the
appendix of the earlier paper \cite{Haj:str} and include basic examples
such as $q$-spheres and $q$-planes.

\subsection*{Acknowledgements} I would like to thank Piotr Hajac and
Tomasz Brzezinski for helpful discussions.

\section{Riemannian geometry with respect to a frame resolution}

Our reformulation of Riemannian geometry is based in the following
lemma about associated bundles. It can be considered as implicit in
\cite{KobNom:fou}, though perhaps not stated explicitly. Let $M$ be
a manifold and $C(M)$ functions on $M$. We work in a smooth setting
throughout the section. Let $\pi:P\to M$ be a principal bundle over
$M$ with structure group $G$, and let $V$ be a left $G$-module. We
recall that there is an associated vector bundle $\CE =P\times_GV$
consisting of equivalence classes in $P\times V$ under the relation
$(p.g,v)\sim (p,g.v)$ for all $p\in P,v\in V$ and $g\in G$. The
projection to $M$ is $\pi(p,v)=\pi(p)$ and the fiber over $x\in M$
is isomorphic to $V$. We recall that sections of $\CE $ are maps
$M\to \CE $ such that following with the projection to $M$ gives
the identity, and that sections may be identified with
pseudotensorial (i.e. $G$-equivariant) functions on $P$,
\[ \Gamma(\CE )\isom C_G(P,V).\]
We extend this now to vectors and forms. We denote by
$\Omega^{-1}(M)$ the space of vector fields on $M$ and, later on,
by $\Omega^1(M)$ the space of 1-forms on $M$. These are sections of
the tangent bundle $TM$ and the cotangent bundle $T^*M$
respectively. A differential from on $P$ is called {\em tensorial} if
it is equivariant and `horizontal' in the sense that it vanishes on the
vertical vector fields
corresponding to the action of $G$ on $P$. Finally, we recall that
a bundle map between bundles over $M$ means a map between the total
spaces forming a commutative triangle with the projections to $M$.

\begin{lemma} The space  $\Omega^1_{\rm tensorial}(P,V)$ of $V$-valued
tensorial 1-forms is in  correspondence with the space of bundle
maps $TM\to \CE $. These in turn are in correspondence with
$C(M)$-module maps
\[ \Omega^{-1}(M)\to C_G(P,V) \]
\end{lemma}
\proof The fiber-wise formulae are as in \cite{KobNom:fou}: given $\theta$
which is $G$-equivariant and horizontal, we define a map
$\tilde\theta:T_xM\to \CE _x$ by
$\tilde\theta(X_x)=[(p,\theta(\tilde X))]$ where $p\in P$ is such
that $\pi(p)=x$ and $\tilde{X}$ is any vector field on $P$ with
horizontal projection $\pi_*(\tilde X)=X$ at $x$ (i.e. any local
lift). One may then check that this construction extends smoothly:
in each open set $U$ choose a local section $\sigma:U\to P$ to
specify the lifts. \endproof

An example is provided by the frame bundle. We recall that a frame
at $x\in M$ is a linear isomorphism $\R^n\to T_xM$ (i.e. a choice
of basis of $T_xM$). The frame bundle $P=FM$ is a principal bundle
over $M$ with structure group $G=GL_n$ and fiber over $x$ given by
the set of all frames at $x$. In a local patch where the bundle is
trivial, we denote by $\pi_2(p):\R^n\to T_{\pi(p)}M$ the
corresponding frame. The canonical 1-form
$\theta\in\Omega^1_{{\rm tensorial}}(P,\R^n)$ is locally
defined by
\[ \theta_p(X)=\pi_2(p)^{-1}\pi_*(X)(\pi(p))\]
for all $X\in TP$. One may check that it is tensorial
and globally defined. Here $V=\R^n$ as a $GL_n$-module and the map
in Lemma~2.1 in this case is an isomorphism $TM\isom
\CE =FM\times_{GL_n}\R^n$. Thus Lemma~2.1 tells us that the tangent
bundle is an associated vector bundle to the frame bundle.

Moreover, the fact that $\theta$ characterises the frame bundle as
observed in \cite{Haj:str}, is now recovered as an immediate
corollary of Lemma~2.1:

\begin{corol} If $(P,\theta), (P',\theta')$ are two $GL_n$-bundles
over $M$ equipped with equivariant tensorial 1-forms with values in
$\R^n$ and such that their induced maps in Lemma~2.1 are
isomorphisms, then $P\isom P'$ and $\theta'$ may be identified with
$\theta$.
\end{corol}
\proof Clearly $P\times_{GL_n}\R^n\isom TM\isom
P'\times_{GL_n}\R^n$ as associated $GL_n$ bundles. Since the
actions on $V=\R^n$ in both cases are the same fundamental
representation of $GL_n$, which is faithful, the structure
constants of $P$ and $P'$ are equivalent, i.e. they are isomorphic
principal bundles.
\endproof

Taking $P=FM$ and $\theta$ the canonical form, any other
$P',\theta'$ with the same structure group $GL_n$ is therefore
isomorphic. But one {\em could} have different examples with
different structure groups. For example, one may take the bundle of
affine frames with structure group $\R^n\lcross GL_n$ or, for a
manifold admitting a metric, the bundle of orthogonal frames with
structure group $O_n$. We call any $(P,G,V,\theta)$ inducing an
isomorphism via Lemma~2.1 a {\em frame resolution} of the tangent
bundle. We now extend these ideas further, to include differential
forms. As a dualisation of Lemma~2.1, we have cf\cite[Sec.
11.14]{KMS:nat}.

\begin{propos} Let $V$ in the setting of Lemma~2.1 be finite-dimensional
and $V^*$ its dual as a $G$-module. $\theta\in\Omega^1_{{\rm
tensorial}}(P,V)$ are in correspondence with bundle maps $\CE ^*\to
T^*M$ and, at the level of sections, with $C(M)$-module maps
\[ C_G(P,V^*)\to \Omega^1(M).\]
\end{propos}
\proof We dualise Lemma~2.1 in a straightforward manner (fiberwise).
Here $\CE ^*=P\times_G V^*$ is the fiber-wise dual. Explicitly, a
$V^*$-valued equivariant function $\phi$ maps to the one form which
at $x$ has values $(\pi^*)^{-1}\phi_p\cdot\theta_p$ where we chose
any $p\in\pi^{-1}(x)$ and identify $\theta_p$ as in the image of
$\pi^*:T_xM\to T_pP$. Here $\cdot$ denotes the evaluation of $V^*$
with $V$. A similar correspondence has been pointed out to us in
\cite[Sec. 11.14]{KMS:nat}; for our purposes we need the
correspondence quite explicitly in the form just given.
\endproof

In the case of the frame bundle, this is an isomorphism and
expresses the cotangent bundle as an associated vector bundle
$T^*M\isom FM\times_{GL_n}\R^{n*}$. Similarly for any frame
resolution $(P,G,V,\theta)$ we see that both 1-forms and vector
fields on $M$ may then be expressed as $V$ or $V^*$-valued
equivariant functions on $P$. In a similar manner, one has in this
case
\eqn{isom1V*}{ \Omega^1(M)\tens_{C(M)}\Omega^1(M)\isom
\Omega^1_{{\rm tensorial}}(P,V^*)}
\eqn{isom1V}{ \Omega^{1}(M)\tens_{C(M)}\Omega^{-1}(M)\isom
\Omega^1_{{\rm tensorial}}(P,V)}
\eqn{isom2V*}{ \Omega^2(M)\tens_{C(M)}\Omega^{1}(M)
=\Omega^2_{{\rm tensorial}}(P,V^*)}
\eqn{isom2V}{ \Omega^2(M)\tens_{C(M)}\Omega^{-1}(M)
=\Omega^2_{{\rm tensorial}}(P,V)}
etc. The canonical form $\theta\in
\Omega^1_{{\rm tensorial}}(P,V)$ corresponds to the constant section of
$\Omega^{-1}(M)\tens_{C(M)}\Omega^1(M)$ given over each point $x$
by the canonical element of $T_xM\tens T^*_xM$.

Now, when $P$ is equipped with a connection, it induces a
connection on $\CE $ and hence a corresponding covariant derivative
$D$ on sections of $\CE $. In the present case it means a covariant
derivative on vector fields, 1-forms etc when these are viewed via
the above isomorphisms as sections of suitable associated vector
bundles. In this way, one obtains the usual formulae of Riemannian
geometry but now as a gauge theory on {\em any} frame resolution.
As such, many of the formulae are in fact more natural.

\begin{lemma} Let $(P,G,V,\theta)$ be a frame resolution of a manifold $M$,
and $\omega$ a connection on $P$. Define the induced
$\nabla:\Omega^1(M)\to\Omega^1(M)\tens_M\Omega^1(M)$ as the
covariant derivative $D:C_G(P,V^*)\to\Omega^1_{{\rm
tensorial}}(P,V^*)$ viewed under the above isomorphisms. Then
$\nabla$ is a derivation with respect to multiplication by
functions on $M$. Moreover,
\[ \pi^*\nabla_X f=\pi^* L_Xf-\phi\cdot \CL_{\tilde X}\theta;\quad
\CL=L+\omega\]
for $\phi\in C_G(P,V^*)$ corresponding to $f\in \Omega^1(M)$ and
$\tilde X$ any lift of a vector field $X$ on $M$. Here $L$ denotes
the Lie derivative and $\cdot$ the evaluation of $V^*$ with $V$.
\end{lemma}
\proof Here $\pi^*f=\phi\cdot\theta$ as in Proposition~2.3, and similarly
$\pi^*\nabla_Xf=(D_{\tilde X}\phi)\cdot\theta= (\tilde{X}(
\phi)+\omega_{\tilde X} \phi)\cdot\theta$ by the similar
isomorphism (\ref{isom1V*}). Note that the covariant derivative on
$C_G(P,V^*)$ is defined by
\[ D\phi=(\id-\Pi_\omega)\extd\phi=D\phi-\tilde\omega(\phi)
=\extd\phi+\omega\phi\]
since $\phi$ is equivariant in the sense $\tilde\xi(\phi)=-\xi\phi$
for all $\xi\in \cg$, the Lie algebra of $G$. Here $\tilde{\xi}$ is
the vector field on $P$ induced by the action of $\xi$ and
$\Pi_\omega$ is the projection on $\Omega^1(P)$ corresponding to
the connection form $\omega$. These steps are the standard
definition of $D\phi$, as in \cite{KobNom:fou}(we recall them
explicitly since we use the quantum group version in Section~3). We
then evaluate against any lift $\tilde X$ of $X$ to a vector field
on $P$. Here $D_{\tilde X}\phi$ is manifestly independent of the
choice of lift since $\omega_{\tilde \xi}=\xi$ for any connection
form, and $\phi$ is equivariant. Also, if $g\in C(M)$ then $gf$
corresponds to $(\pi^*g)\phi$. Hence $\pi^*\nabla
gf=(D((\pi^*g)\phi))\cdot\theta=\pi*^(\extd g)\phi\cdot
\theta+(\pi^*g)(D\phi)\cdot\theta$, i.e. $\nabla$ is a derivation.
Finally, writing $\tilde X(\phi)=L_{\tilde X}\phi$ and the Leibniz
property of the Lie derivative (and moving $\omega$ to act on $V$
rather than $V^*$ by $(\xi\phi)\cdot\theta=-\phi\cdot\xi\theta$)
gives the alternative expression for $\nabla_X f$ as stated.
\endproof

We call the operation $\CL=L+\omega$ on form-sections  the {\em
covariant Lie derivative}. We see that $\CL\theta$ measures the
deviation of the covariant derivative from the Lie derivative on
forms.  Finally, we define the covariant derivative on vector
fields by
\eqn{nabvec}{ \<\nabla_XY,f\>=X(\<Y,f\>)-\<Y,\nabla_Xf\>}
for all vector fields $X,Y$ on $M$ and all $f\in\Omega^1(M)$.

\begin{propos} Define the torsion tensor $T$ as the (1,2)-tensor
on $M$ corresponding
to $D\theta\in\Omega^2_{{\rm tensorial}}(P,V)$ under the above
isomorphisms. Then
\[ \nabla\wedge f=\extd f -\<T,f\>\]
for all $f\in\Omega^1(M)$, which is equivalent to the usual
definition of the torsion tensor.
\end{propos}
\proof We first show that the equation shown for $T$ is equivalent
to the usual definition
\[ \nabla_X Y-\nabla_Y X=[X,Y]+T(X,Y)\]
for all vector fields $X,Y$ on $M$. Using (\ref{nabvec}), this is
equivalent to
\[X(\<Y,f\>)-\<Y,\nabla_Xf\>-Y(\<X,f\>)+\<X,\nabla_Yf\>-\<[X,Y],f\>
=\<T(X,Y),f\>.\]
On the other hand, using the notation $i_X$ for interior product
with forms, $L_X$ for Lie derivative, and the identities
$L_X=i_X\extd+\extd i_X$, $[L_X,i_Y]=i_{[X,Y]}$, we have
\eqn{iyix}{i_Yi_X\extd f=-i_Y\extd i_X f+i_Y L_Xf=
X(\<Y,f\>) -Y(\<X,f\>)-\<[X,Y],f\>} as required. Here
$i_Yi_X(\nabla\wedge f)=(\nabla
\wedge f)(X,Y)=\<Y,\nabla_Xf\>-\<X,\nabla_Yf\>$

We now let $T$ correspond to $D\theta$ by (\ref{isom2V}).
Equivalently, $T(X,Y)$ is a vector field corresponding under the
induced isomorphism in Lemma~2.1 to $(D\theta)_{\tilde X,\tilde
Y}\in C_G(P,V)$. Hence $\pi^*\<T(X,Y),f\>=\phi\cdot
(D\theta)_{\tilde X,\tilde Y}=i_Yi_X \phi\cdot D\theta$.

Once these isomorphisms are understood, the computation in the
tensorial form language is immediate:  $\pi^*\nabla\wedge
f=(D\phi)\wedge\theta=(\extd
\phi)\wedge\theta+\omega\phi\wedge\theta= \extd(\phi\cdot\theta)
-\phi\cdot \extd\theta-\phi\cdot\omega\wedge\theta=\pi^*\extd f
-\phi\cdot D\theta$
as required. \endproof

The operation $\nabla\wedge$ is the {\em covariant exterior
derivative} and we see that the torsion $T$ measures its difference
from the usual exterior derivative on 1-forms. We can similarly
treat the Riemannian curvature in terms of differential forms and
sections.

\begin{propos} Define the curvature tensor of $\nabla$ by $Rf$ the 1-form
corresponding to $F\phi$ where $\Omega\in\Omega^2_{\rm
tensorial}(P,\cg)$ is the curvature of the connection form $\omega$
and $f\in\Omega^1(M)$ corresponds to $\phi$. Then
\[R(X,Y)f=[\nabla_X,\nabla_Y]f-\nabla_{[X,Y]}f\]
for all vector fields $X,Y$ on $M$.
\end{propos}
\proof It is easy to check that the stated formula for $R(X,Y)$
acting on forms
is equivalent to its usual definition as an operator on vector
fields, when the two are related by $\<Z,R(X,Y)f\>=-\<R(X,Y)Z,f\>$
for all $X,Y,Z$ and $f$. We use (\ref{nabvec}) repeatedly to
establish this. Next, from the definition of $\nabla$ above,
$\pi^*(\nabla_X\nabla_Yf-\nabla_Y\nabla_Xf-\nabla_{[\tilde X,\tilde
Y]}f)= (([D_{\tilde X},D_{\tilde Y}]-D_{[\tilde X,\tilde
Y]})\phi)\cdot\theta$. To compute this, note that
\align{\tilde
X(\omega_{\tilde Y}\phi)\equad &&=L_{\tilde X}i_{\tilde
Y}(\omega\phi)=i_{\tilde Y}L_{\tilde X}(\omega\phi)+i_{[\tilde
X,\tilde Y]}\omega\phi\\ &&=\omega_{\tilde Y}\tilde
X(\phi)+i_{\tilde Y}(L_{\tilde X}\omega)\phi+\omega_{[\tilde
X,\tilde Y]}\phi=\omega_{\tilde Y}\tilde X(\phi)+i_{\tilde
Y}i_{\tilde X}(\extd \omega)\phi+i_{\tilde Y}(\extd i_{\tilde
X}\omega)\phi+\omega_{[\tilde X,\tilde Y]}\phi\\ &&=\omega_{\tilde
Y}\tilde X(\phi)+i_{\tilde Y}i_{\tilde X}(\extd \omega)\phi+\tilde
Y(\omega_{\tilde X}\phi)-\omega_{\tilde X}\tilde
Y(\phi)+\omega_{[\tilde X,\tilde Y]}\phi.} We used the usual
formulae for $[L,i]$ and $L=i\extd+\extd i$, this time in $P$.
Hence
\align{[D_{\tilde X},D_{\tilde Y}]\phi\equad &&=D_{\tilde X}(\tilde Y(\phi)
+\omega_{\tilde Y}\phi)-(X\swap Y)=\tilde X(\tilde
Y(\phi))+\omega_{\tilde X}\tilde Y(\phi)+\tilde X(\omega_{\tilde
Y}\phi)+\omega_{\tilde X}\omega_{\tilde Y}\phi-(X\swap Y)\\
&&=[\tilde X,\tilde Y]\phi+\omega_{[\tilde X,\tilde Y]}\phi+
i_{\tilde Y}i_{\tilde X}(\extd
\omega+\omega\wedge\omega)\phi=\Omega_{\tilde X,\tilde Y}\phi} as required.
\endproof

One may also define the exterior covariant derivative on 2-forms as
corresponding to $D$ on $\Omega^1(P,V^*)$, and then
$\nabla\wedge\nabla\wedge f=R\wedge f$ holds, cf \cite{KMS:nat}.
Moreover, to complete the picture, we can also consider the metric
under
\eqn{isomV*V*}{\Omega^1(M)\tens_{C(M)}\Omega^1(M)
\isom C_G(P,V^*)\tens_{C(M)}C_G(P,V^*)\isom C_G(P,V^*\tens V^*)}
as an equivariant function on $P$ with values in $V^*\tens V^*$.
The first isomorphism here is that in Proposition~2.3 applied to
each tensor factor over $C(M)$. The second is given by pointwise
product in $P$ and is clearly an isomorphism for trivial bundles
(where $C_G(P,V^*)=C(M,V^*)$ etc.), and hence holds generally by
local triviality of a general bundle.

\begin{propos} Let $\nabla_X$ on any $g\in\Omega^1(M)
\tens_{C(M)}\Omega^1(M)$ be defined
as corresponding to $D:C_G(P,V^*\tens V^*)
\to\Omega^2_{\rm tensorial}(P,V^*\tens V^*)$ under the above isomorphisms.
Then $\nabla_X$
is the extension of $\nabla_X$ on $\Omega^1(M)$ to the tensor
product as a derivation.
\end{propos}
\proof We make a similar computation to that in Lemma~2.4. Thus,
if $g\in \Omega^1(M)\tens_{C(M)}\Omega^1(M)$ is viewed as $\eta\in
C_G(P,V^*\tens V^*)$ then
\align{\pi^*\nabla_X g\equad &&=(D_{\tilde X}\eta)\cdot(\theta\tens\theta)
=\pi^*L_Xg-\eta\cdot\CL_{\tilde X}(\theta\tens\theta)\\
&&=\pi^*((L_Xg_\alpha) g^\alpha+g_\alpha L_X
g^\alpha)-\eta\cdot(\CL_{\tilde
X}\theta\tens\theta+\theta\tens\CL_{\tilde X}\theta),}where
$g=g_\alpha\tens_{C(M)}g^\alpha$, say (summation understood). We
used the derivation property of the Lie derivative and the fact
that $\cg$ (where lies the output of $\omega$) acts like a
derivation on the tensor product representation $V^*\tens V^*$.
Finally, by the alternative expression for $\nabla$ on
$\Omega^1(M)$ in Lemma~2.4, we obtain
\eqn{gder}{ \nabla_Xg=(\nabla_Xg_\alpha)g^\alpha+g_\alpha\nabla_X g^\alpha}
as stated. Note that evaluation against $Y,Z$ and (\ref{nabvec})
shows that (\ref{gder}) is equivalent to the more conventional
definition of $\nabla_X g$ as\cite{KobNom:fou}
\[ (\nabla_X g)(Y,Z)=X(g(Y,Z))-g(\nabla_XY,Z)-g(Y,\nabla_XZ). \]
\endproof

Finally, in the case where the frame resolution is a trivial
bundle, the form $\theta$ itself corresponds to the `abstract
$n$-bein' or {\em $V$-bein} $\vb\in\Omega^1(M,V)$ via
$\theta(p)=\pi_2(g)^{-1}\pi^*\vb(p)$. The space $C_G(P,V^*)\isom
C(M)\tens V$ by $\phi(p)=\pi_2(p)\pi^*\psi$ where $\psi\in
C(M,V^*)$ is a `matter field' with values in $V^*$. These
particular `matter fields' correspond to 1-forms $f\in
\Omega^1(M)$ by $f=\psi\cdot \vb$. Similarly, a
metric corresponds to `matter field' $\eta\in C(M,V^*\tens V^*)$ by
$g=\eta\cdot (\vb\tens \vb)$. Moreover, as usual, a connection
$\omega$ corresponds to gauge field $A\in\Omega^1(M,\cg)$ by
\[\omega(p)=\pi_2(p)^{-1}\pi^*A(p)\pi_2(p)+\pi_2(p)^{-1}\extd\pi_2(p)\]
with the usual abuses of notation for the second term. This is how
the quantities above look in terms of the usual `matter fields' and
`gauge potentials'. Most manifolds do not admit trivial frame
resolutions, but these are also the local formulae for each patch
of a nontrivial frame resolution. One has similar formulae to the
general case, for example
\eqn{locnab}{ \nabla_Xf=L_Xf-\psi\cdot\CL_Xh,\quad \CL=L+A}
if 1-form $f$ corresponds to matter-field $\psi$. Here $\CL$ can be
called the {\em covariant Lie derivative} on $M$. In particular,
the Levi-Civita connection corresponds to (an example of) a gauge
field leaving covariantly constant a matter field $\eta$ with
values in $V^*\tens V^*$ and a 1-form $\vb$ with values in $V$.

This completes our formulation of the main properties of the
covariant derivative on forms in terms of connections on an
arbitrary frame resolution. In the general case we could call a
connection `Levi-Civita' {\em with respect to a framing and a
metric} if $D\theta=0$ (torsion free) and $D\eta=0$ (metric
compatible). We should not expect existence and uniqueness however,
as this depends strongly on a particular choice of resolution.
Instead, one should reverse the logic and think of the frame
resolution as part of the input data and already playing in part
the role of choice of a metric.  The choice of this and a
connection on it specifies a covariant derivative {\em in lieu} of
a metric. In the standard $O_n$ frame bundle case we know that
connections correspond to metrics in the usual way, while a
subgroup would only allow a subset of connections or a subset of
corresponding metrics. Under this kind of correspondence, more
general $G$ are also possible, corresponding to locally preserving
different kinds of `generalised metrics'. For example, we may take
$G$ symplectic and symplectic forms in the role of metric. We could
also consider manifolds equipped with, say, $E_6$ frame resolution
or with resolution by infinite-dimensional groups, even without
consideration of any metric.

Moreover, for any fixed $G$ we can consider frame resolutions of
$\Omega^1(M)$ via different representations $V$. For example, the
spinor representation in the $O_4$ case leads ultimately to the
spin bundle and Dirac operator. As a more novel application of this
idea we consider now  the representation $V^*$ conjugate to any
$V$. We show that this leads to a natural `self-dual' formulation
and slight generalisation of Riemannian geometry. We reformulate a
(not necessarily symmetric) metric $g$ as corresponding under
(\ref{isom1V*}) to $\gamma\in\Omega^1_{\rm tensorial}(P,V^*)$ and
develop the theory symmetrically between $\gamma$ and $\theta$ as
defining conjugate frame resolutions, namely associated to $V,V^*$
respectively.

\begin{corol} Given a frame resolution $(P,G,V,\theta)$, a 2-cotensor $g$ is
nondegenerate as a map $g:\Omega^{-1}(M)\isom\Omega^1(M)$ {\em iff}
the corresponding $\gamma$ makes $(P,G,V^*,\gamma)$ into another
frame resolution, called the {\em dual frame resolution}.
\end{corol}
\proof Under the  isomorphism $s_\theta$ of Corollary~2.2, an
isomorphism $g:\Omega^{-1}(M)\to
\Omega^1(M)$ (by $X\mapsto g(X,\ )$) is equivalent to an
isomorphism $\Omega^{-1}(M)\to C_G(P,V^*)$. It sends $X$ to
$\gamma_{\tilde X}$ (for any lift $\tilde X$) since $\gamma_{\tilde
X}\cdot\theta=\pi^*g(X,\ )$. By Lemma~2.1, this is equivalent to a
frame resolution $(P,G,V^*,\theta')$ where $\theta'$ is a tensorial
$V^*$-valued 1-form on $P$ such that $\Omega^{-1}(M)\to C_G(P,V^*)$
is given by mapping $X$ to $\<\theta',\tilde X\>$, i.e.
$\theta'=\gamma$.

Note that, as a corollary, we have an explicit correspondence
\eqn{phif}{f\in\Omega^1(M)\quad \swap\quad \phi=\gamma_{\tilde{g^{-1}(f)}}}
in Proposition~2.4 etc., giving explicit formulae for the
expressions there.
\endproof

We can also reverse the roles of $\theta$ and $\gamma$ (and $V^*$
and $V$) in all of the above, regarding $\gamma$ as the frame
resolution and $\theta$ as corresponding to a generalised metric.
From this point of view it is natural to replace $D\theta=0$ and
$\nabla g=0$ by more symmetric `self-dual' conditions $D\theta=0$,
$D\gamma=0$.

\begin{propos} Let $(P,G,V,\theta)$ be a frame resolution and $g$
a nondegenerate
2-cotensor as `generalised metric', viewed as a corresponding dual
resolution $(P,G,V^*,\gamma)$. We define the {\em cotorsion form}
$\Gamma\in
\Omega^2(M)\tens_{C(M)}\Omega^1(M)$ by
\[\Gamma(X,Y)=(\nabla_X g)(Y,\ )-(\nabla_Y g)(X,\ )+g(T(X,Y))\]
for vector fields $X,Y$ on $M$. Then (i) $\Gamma=g(T_\gamma)$ where
$T_\gamma$ is the torsion of $\gamma$ in the dual frame resolution,
and (ii) $\Gamma$ corresponds under (\ref{isom2V*})
 to $D\gamma$.
\end{propos}
\proof  It is easy to see from the preceding corollary that $D\gamma$
corresponds
under (\ref{isom2V*}) to $g(T_\gamma)$. We prove part (ii). Using
the identity (\ref{iyix}) applied on $P$ (and applied to lifts
$\tilde X,\tilde Y$), we have
\align{\<Z,\Gamma(X,Y)\>\equad&&=(i_{\tilde Y}i_{\tilde X} D\gamma)
\cdot \theta_{\tilde Z}\\
&&=\tilde X(\gamma_{\tilde Y})\cdot\theta_{\tilde Z}-\tilde
Y(\gamma_{\tilde X})\cdot\theta_{\tilde Z}-\gamma_{[\tilde X,\tilde
Y]}\cdot\theta_{\tilde Z}+(\omega_{\tilde X}\gamma_{\tilde
Y})\cdot\theta_{\tilde Z}-(\omega_{\tilde Y}\gamma_{\tilde
X})\cdot\theta_{\tilde Z}\\ &&= \<\nabla_X g_Y,Z\>-\<\nabla_Y
g_X,Z\>-g([X,Y],Z)\\ &&=(\nabla_X g)(Y,Z)-(\nabla_Y
g)(X,Z)+g(T(X,Y),Z).} Here $g_Y=g(Y, \ )\in \Omega^1(M)$
corresponds to $\gamma_{\tilde Y}$ as in the preceding corollary,
and its covariant derivative by $\nabla_X$ therefore corresponds to
$D_{\tilde X}\gamma_{\tilde Y}$. The last line comes from
$\<\nabla_X g_Y,Z\>=X(g(Y,Z))-g(Y,\nabla_XZ)$ as in (\ref{nabvec}),
the usual expression for the torsion tensor as in Proposition~2.5,
and the usual expression for $\nabla_X g$ as in Proposition~2.7.
\endproof

We note, using Proposition~2.5 and the Leibniz property as in
Proposition~2.7, that we may also write the cotorsion form as
\eqn{cotor}{g(T_\gamma)= \nabla\wedge g_\alpha\tens g^\alpha+g(T)
-g_\alpha\wedge\nabla \tens g^\alpha
=\extd g_\alpha\tens g^\alpha-g_\alpha\wedge \nabla\tens g^\alpha.}
This gives an immediate corollary:

\begin{corol} In the setting of Proposition~2.9, suppose that the torsion vanishes.
Then the cotorsion vanishes {\em iff} $\extd g=0$, where $\extd$ is
extended to 2-cotensors as a graded-derivation.
\end{corol}
\proof  We project (\ref{cotor}) to $\Omega^3(M)$. Then the projection
of $g(T_\gamma)$ is $\extd g+g_\alpha \wedge\<T,g^\alpha\>$.
\endproof

If we call $g(\ ,T)\in\Omega^1(M)\tens_{C(M)}\Omega^2(M)$ the {\em
cotorsion form}, the corollary says that the skew-symmetrized
cotorsion form and torsion form differ by $\extd g$. For example,
if $g$ is antisymmetric, $\extd g$ is the De Rahm exterior
differential  the vanishing of cotorsion is the same as saying that
$g$ is a symplectic 2-form. I.e. symplectic geometry is naturally
included in our generalisation. Finally, when the frame resolution
bundle is trivial, $\gamma$ corresponds to $\cvb\in
\Omega^1(M,V^*)$ (which we call a {\em $V$-cobein}), and $g=\cvb\tens\vb$.
The vanishing of torsion and cotorsion
with respect to a gauge field $A$ is then the symmetrical condition
\eqn{0torcotor}{ D_A\vb=0,\quad D_A\cvb=0.}

As an example, if $G$ is any semisimple Lie group, we consider
$G\times G$ as a principal $G$-bundle (where $G$ acts on the right
factor from the right). We take $V=\cg$ with the adjoint action. On
$G$ is a canonical Maurer-Cartan form with values in $\cg$. We take
this for $\vb$ and we let $\cvb=\eta\circ\vb$ where $\eta$ is the
Killing form. This defines the usual metric on $G$ (Riemannian in
the compact case). One also has a natural gauge field $A$ inducing
a covariant derivative with vanishing torsion and cotorsion. It is
given by $A$ equal to $\h$ the Maurer-Cartan form. We interpret
$\extd A+A\wedge A=0$ variously as (\ref{0torcotor}) and zero
curvature.  We can similarly treat homogeneous spaces $G/H$.

Finally, as well as suggesting such natural generalisations of
conventional Riemannian geometry, the gauge field formulation also
suggests more radical approaches to quantisation of the graviton.
For example, one may regard the frame resolution data
$(P,G,V,\theta,\gamma)$ as part of the `background manifold' data
in which our particles move, considering only the connection or
gauge field $A$ part as dynamic. There are clearly a number of
interesting directions suggested by the above formulation, some of
which will be explored elsewhere. The main point for our present
purposes, however, is that it provides a clean `co-ordinate free'
and `differential form' approach suitable for generalisation to
quantum group non-commutative geometry, to which we now turn.

\section{Quantum group Riemannian geometry}

Motivated by the above formulation of classical Riemannian
geometry, we now let $M$ be a possibly non-commutative unital
algebra over a general ground field $k$. We use the theory of
quantum principal bundles $P$ with quantum structure group $\qg $,
as introduced in \cite{BrzMa:gau}. Here $\qg $ coacts by $\Delta_R$
on $P$, $M=P^\qg $ is the fixed point subalgebra, $P$ is assumed to
flat as an $M$-module and the extension is assumed to be
Hopf-Galois. These definitions are somewhat like the classical ones
but with arrows reversed since our spaces $M$ etc are replaced by
algebras playing the role of their ring of functions. We use the
universal differential calculi $\Omega^1M,\Omega^1P$ etc associated
to any unital algebra. Here $\Omega^1M\subset M\tens M$ is the
kernel of the product map and $\extd:M\to \Omega^1M$ is $\extd
m=1\tens m-m\tens 1$. Further, one finds\cite{BrzMa:gau} effectively that
the role of local triviality in the theory of connections can be played 
by the assumption that the map $\chi:P\tens_MP\to 
P\tens \qg $ defined by
the descent to $P\tens_MP$ of
$\tilde\chi=(\cdot\tens\id)(\id\tens\Delta_R)$ is invertible. This
`Hopf-Galois' condition is familiar in Hopf algebra theory and has been
considered as a natural `topological' requirement \cite{Sch:pri} even
without the differential calculus and gauge theory for this setting
introduced in \cite{BrzMa:gau}. We assume that the quantum group $\qg$ has invertible antipode.

As shown in \cite{Ma:war95} and cf\cite{BrzMa:gau}\cite{Brz:tra},
if $V$ is a right $\qg $-comodule with invariant `unit' element
$1\in V$, we have an associated fiber bundle $\CE =(P\tens V)^\qg $
containing $M$ as $M\tens 1$. The sections of $\CE $ are the unital
left $M$-module maps $\CE \to M$ and are in 1-1 correspondence with
to unital equivariant maps $V\to P$. To this theory we now add a
description of $\Omega^nM$ in terms of $\CE $. These kind of
results do not actually need a unit element in $V$ and we drop this
for the purposes of the present paper. We also recall from
\cite{BrzMa:gau} that $\theta:V\to\Omega^1P$ is right strongly
tensorial (r.s.t.) if it is equivariant and its image lies in
$P(\Omega^1M)$.

\begin{lemma} Right strongly tensorial $\theta:V\to P\Omega^nM$ are in
1-1 correspondence with maps left $M$-module maps $\CE \to
\Omega^nM$.
\end{lemma}
\proof In principle, we identify the sections of the bundle $\CE ^*
=(P\tens V^*)^\qg $
with equivariant maps $V^*\to P$, i.e. with elements of $\CE
=P\tens V$. Here $V^*$ is the right dual, i.e. has a coaction such
that the evaluation map $V\tens V^*\to \C$ is an intertwiner. In
practice, we now proceed directly with $E$ and without assuming
$V^*$. Given $\theta$, we define $\tilde s_\theta:P\tens V\to
P\tens M^{\tens n}$ by $\tilde s_\theta(p\tens v)=p\theta(v)$ and
verify that its restriction $s_\theta:\CE
\to\Omega^nM$ to $\CE
=(P\tens V)^\qg $ indeed has its image in $\Omega^nM$. To see this,
apply the coaction $\Delta_R$ to $P$ in $P\tens M^{\tens n}$.
Conversely, given $s:\CE
\to \Omega^nM$, let $\theta_s(v)=\chi^{-1}_\alpha(1\tens
S^{-1}v\bt).s(\chi^{-1\alpha}(1\tens S^{-1}v\bt)\tens v\bo)$ which
one may verify is well-defined and lies in $P\Omega^nM$. Here
$\Delta_R(v)=v\bo\tens v\bt$ (summation understood) is the coaction
on $v$ and $\chi^{-1}=\chi^{-1}_\alpha\tens\chi^{-1\alpha}$ is
another notation (summation over terms labeled by $\alpha$
understood). The proof follows exactly the same lines as the $n=0$
case in \cite{Brz:tra}\cite{Ma:war95}.
\endproof

Motivated by Corollary~2.2 we define a {\em frame resolution} of
$\Omega^1M$ as follows:

\begin{defin} A frame resolution of $(M,\Omega^1M)$ is a quantum group
principal bundle $P(M,\qg )$ over $M$, a unital right $\qg
$-comodule $V$ and a right strongly tensorial form $\theta:V\to
P\Omega^1 M$ such that the map $s_\theta:\CE \to \Omega^1M$ is an
isomorphism.
\end{defin}

Note that we do not fix $V$ or $\qg $ here, so this is not
necessarily the frame bundle in the classical case.  However, the
moral of our formulation of classical Riemannian geometry in
Section~2 is precisely that it does not really matter which frame
resolution we take as all of them  achieve the desired result that
the (co)tangent bundle is expressed as an associated vector bundle.
We can go on to define the `frame bundle' as a particular frame
resolution (e.g. in some sense minimal such that there is a unique
torsion free connection for every metric), but we do not need to do
this (it will be attempted elsewhere).

We can now proceed to conclude further isomorphisms
\eqn{qisom1V}{ \id\tens s_\theta:((\Omega^1M)P\tens V)^\qg
\isom \Omega^1M\tens_M\Omega^1M}
(i.e. with left strongly tensorial forms $V^*\to (\Omega^1M)P$),
etc.

\begin{propos} Let $\omega$ be a left strong connection on $P$ in
the sense of preserving left strong tensoriality. Then the
covariant derivative $D=(\id-\Pi_\omega)\extd$ on the associated
bundle $\CE ^*$ corresponds to a map $\nabla:\Omega^1M\to
\Omega^1M\tens_M\Omega^1M=\Omega^2M$ obeying the derivation property
with respect to left-multiplication by $M$. Explicitly,
\[ \nabla =1\tens\id-(\id\tens\theta)\circ s_\theta^{-1}
-s^{-1}_{\theta \alpha}\bo\cdot\omega(s^{-1}_{\theta \alpha}\bt)\cdot
\theta(s^{-1 \alpha}_{\theta})\]
where $s^{-1}_{\theta \alpha}\tens s^{-1\alpha}_\theta$ (summation
over terms labeled by $\alpha$ understood) denotes the output of
$s_\theta^{-1}$ after acting on $\Omega^1M$.
\end{propos}
\proof By definition, a connection on $P$ is called (left) strong if it sends
pseudotensorial forms $f:V^*\to P$ (which are automatically left
and right tensorial) to left strongly tensorial forms $Df:V^*\to
(\Omega^1M)P$. (The need for such a restriction is explained in
\cite{BrzMa:gau} and it is studied in \cite{Haj:str}). We view $f$
as a section of $\CE ^*$ or an element of $\CE $, and hence via the
above isomorphisms we obtain the corresponding covariant derivative
$\nabla$ on forms as $\nabla=(\id\tens
s_\theta)((\id-\Pi_\omega)\extd
\tens\id)s_\theta^{-1}$. Putting in the form of $\extd$ and
$\Pi_\omega(p\tens p')=pp'\bo\omega(p'\bt)$ for all $p,p'\in P$
from \cite{BrzMa:gau} gives the explicit expression shown for
$\nabla$. From this it is immediate that the derivation property
$\nabla(f.w)=f\nabla w+\extd f\wedge w$ holds for all $f\in M$ and
$w\in\Omega^1M$. Note that the product $\wedge$ in the universal
case is the product of the adjacent copies of $M$.
\endproof

This constructs the covariant derivative on $\Omega^1M$, which is
the starting point of quantum group Riemannian geometry. Next we
define the torsion tensor (motivated by Proposition~2.5) as
$\extd-\nabla$. We have a problem, however, that $\theta$ is
right-strongly tensorial not left-strongly tensorial, so $D\theta$
is not strongly tensorial. This problem is resolved as follows.
Note that assuming $\qg $ has invertible antipode then$\bar
\qg =\qg ^{\rm op}$ (with the opposite product) is also a Hopf algebra
and $\bar P=P$ as an algebra but with the coaction
\eqn{deltaL}{\Delta_L(p)=S^{-1}p\bt\tens p\bo}
makes $\bar P$ into a left comodule algebra under $\bar \qg $.
Throughout the paper, whenever we consider both left coactions and
right coactions $\Delta_Rp=p\bo\tens p\bt$ on the same object, they
will always be related like this. We define
$\bar\Pi_\omega:\Omega^1P\to\Omega^1P$ by the restriction to
$\Omega^1P$ of
\[ \bar\Pi_\omega(p\tens p')=\omega(p\tildebo)p\tildebt p'\]
for $p,p'\in P$ and $\Delta_Lp=p\tildebo\tens p\tildebt$ is a
notation (summation understood). This is exactly a left-handed
version of the formula for $\Pi_\omega$ in \cite{BrzMa:gau}. We
define $\bar D=(\id-\bar\Pi_\omega)\extd$ on forms on $P$.

\begin{propos} The following three are equivalent:

(i)  $\omega$ is a left strong connection on $P$
\eqn{leftstr}{{\it (ii)}\qquad\qquad (\Delta_L\tens\id)\circ\omega
=\id\tens1\tens 1
-1\tens1\tens 1\eps+(\id\tens\omega)\circ\Delta}
\eqn{rightstr}{{\it (iii)}\qquad\qquad(\id\tens\Delta_R)\omega
=1\tens1\tens\id-1\tens1\tens1\eps+(\omega\tens\id)\circ\Delta}
Moreover, in this case $\bar D$ preserves right strongly tensorial
form.
\end{propos}
\proof By definition, $\omega$ is left strong if it preserves left
strongly tensorial
forms. As explained in \cite{Ma:war95}, this is equivalent to
sending   the identity map $P\to P$ regarded as a (left) strongly
tensorial form to a left strongly tensorial form, which is the
condition studied in \cite{Haj:str}. The explicit form of this with
the universal calculus in \cite{Ma:war95}  can be massaged further
(using the comodule algebra property of $\Delta_R$) as
\[  p\bo\omega_\alpha(p\bt\t)\bo\tens p\bt\o\omega_\alpha(p\bt\t)\bt\tens
 \omega^\alpha(p\bt\t) =p\tens 1\tens 1
-\Delta_Rp\tens 1+p\bo\omega_\alpha(p\bt)\tens1\tens\omega^\alpha(p\bt)\]
where  $\omega_\alpha\tens\omega^\alpha=\omega$ is a notation.
Multiplying all expressions here from the left by $p'\in P$, we
note that they can all be factored through the map $\chi(p'\tens_M
p)=p'\Delta_R(p)$, which is invertible by the Hopf-Galois
condition. So, applying to $\chi^{-1}(1\tens h)$ for $h\in \qg $,
we obtain the simpler but equivalent condition
\eqn{leftstra}{\omega_\alpha(h\t)\bo\tens
h\o\omega_\alpha(h\t)\bt\tens\omega^\alpha(h\t)=\eps(h)1\tens1\tens1-1\tens
h\tens 1+\omega_\alpha(h)\tens 1\tens\omega^\alpha(h).} Finally, we
apply this equation to $h\t$ of $Sh\o\tens h\t$ and multiply the
free $Sh\o$ in from the left.  This gives the simpler equation
\[  (\Delta_R\tens\id)\omega(h)
=1\tens   Sh\tens 1-\eps(h)1\tens 1\tens 1
+\omega_\alpha(h\t)\tens Sh\o\tens \omega^\alpha(h\t). \] This is
the condition stated in terms of $\Delta_L$ and is equivalent so
long as $S$ is invertible. On the other hand, $\omega$ as a
connection is $\Ad$ invariant. Using this fact, we may convert the
coaction on the first output of $\omega$ on the left hand side of
(\ref{leftstra}) to a coaction on the second output of $\omega$:
\eqn{rightstra}{ \omega_\alpha(h\o) \tens
\omega^\alpha(h\t)\tildebo\cdot_{\rm op}h\t\tens \omega^\alpha(h\t)\tildebt
=\eps(h)1\tens1\tens1-1\tens
h\tens 1+\omega_\alpha(h)\tens 1\tens\omega^\alpha(h)} as
equivalent to (\ref{leftstra}) by $\Ad$-invariance of $\omega$.
This is the same condition as (\ref{leftstra}) with left-right
reversed and $\qg $ replaced by $\bar \qg $. We can also apply it
to $h\o$ of $h\o\tens h\t$ and multiplying in the $h\t$ gives us
the equivalent expression (\ref{rightstr}).

Moreover, if $\bar P(\bar \qg ,M,\Delta_L)$ is a left handed
quantum principal bundle then we can immediately conclude the
result by left-right symmetry: view a right strongly tensorial form
on $P$ equivalently as a strongly tensorial form on $\bar P$: the
condition that such forms are preserved under $\bar D$ is exactly
the $\Delta_L$ version of (\ref{leftstra}), which have seen is
equivalent to $\omega$ left strong. More generally, without
assuming $\bar P$ is Hopf-Galois we define $\bar D$ directly by the
same formulae as in the Hopf-Galois case and verify that it
preserves right strongly tensorial forms. \endproof

The proposition suggests the definition of a quantum bundle as {\em
bicovariant} if both $P,\bar P$ obey the Hopf-Galois condition. In
this case, the above lemma asserts that a connection is
(left)strong on $P$ {\em iff} it is (right)strong on $\bar P$. We
do not need to assume this at the moment since we are interested
only in working on $P$.

\begin{propos} The torsion tensor $T:\Omega^1M\to\Omega^1M\tens_M\Omega^1M$
defined as $T=\extd-\nabla$ corresponds under the isomorphism in
Lemma~3.1 to the right strongly tensorial form $\bar D\theta:V\to
P\Omega^2M$. Explicitly,
\[ \bar D\theta=1\tens\theta-\theta_\alpha\tens1\tens\theta^\alpha+\theta\tens 1
+\cdot\circ(\omega\tens\theta)\circ\Delta_L\] where
$\theta=\theta_\alpha\tens\theta^\alpha$ (say) and $\cdot$
multiplies the adjacent copies of $P$ between $\omega$ and
$\theta$, and $\Delta_L$ is related as in (\ref{deltaL}) to the
original right coaction on $V$.
\end{propos}
\proof The exterior derivative on $\Omega^1M\subset M\tens M$ is
the restriction
of $\extd(m\tens n)=1\tens m\tens n-m\tens 1\tens n+m\tens n\tens
1$. Subtracting   $\nabla(m\tens n)$ given in Proposition~3.5 and
precomposing all terms with $s_\theta$ gives $(\extd-\nabla)\circ
s_\theta:\CE \to \Omega^2M$. By Lemma~3.1 it has the form
$\cdot\circ(\id\tens Y)$ for some  1-form $Y:V\to P\Omega^2M$.
Computing $Y$ gives the same expression as that stated for $\bar
D\theta$. On the other hand, this is what one obtains as
$(\id-\bar\Pi_\omega)\extd\theta$ by similar (but left-handed)
computations to those in \cite{BrzMa:gau}. \endproof

Clearly, a torsion free connection with respect to a given frame
resolution is one where $\bar D\theta=0$. Let us note that in the
case of the universal calculus where
$\Omega^2M=\Omega^1M\tens_M\Omega^1M$, this condition is much
stronger than it would be classically or with non-universal
calculi.  Since $\nabla$ coincides with $\nabla\wedge$, we see that
the vanishing of torsion implies a unique covariant derivative,
namely $\nabla=\extd$. One may also see by applying $s_\theta^{-1}$
to the explicit formula for $\bar D\theta$ that  $\omega$ composed
with the left-handed coaction $\Delta_L:V\to \qg \tens V$ is fully
determined by $\bar D\theta$. One has uniqueness of the composition
for prescribed torsion. Moreover, by considering $\bar D^2\theta$,
one finds that the curvature of such an $\omega$ must vanish. These
are intrinsic limitations of the universal calculus.

Next, we can  go on and define  metric-compatibility as $\nabla
g=0$, where $\nabla$ is suitably extended to
$g\in\Omega^1M\tens_M\Omega^1M$, e.g. as a derivation if we follow
the line of Proposition~2.7.  On the other hand, Corollary~2.8
suggests a different formulation as follows. We consider $g$
equivalently as $\gamma\in (\Omega^1M)P\tens V)^\qg $ via the
isomorphisms above non-degenerate in the sense of providing a dual
(and left-handed) frame resolution via $V^*$. Thus,
\eqn{ggamma}{ g=\gamma(f^a)\theta(e_a)\in \Omega^1M\tens_M
\Omega^1M\subset M\tens M\tens M,}
where we assume the existence of the canonical element or
coevaluation $f^a\tens e_a=\coev\in V^*\tens V$  for the duality
pairing of $V^*$ with $V$ (here $\{e_a\}$ is a basis of $V$ and
$\{f^a\}$ a dual basis.) When $V$ is infinite-dimensional one can
consider (\ref{ggamma}) formally as a power-series or alternatively
one can replace $\gamma$ by a {\em left-handed} frame resolution
$\theta_L:V\to (\Omega^1M)P$ (i.e. a left strongly tensorial form
such that $s_{\theta_L}(v\tens p)=\theta_L(v)p$ induces an
isomorphism $V\tens M\isom\Omega^1M$) and an invariant element
$\eta=\eta\uo\tens\eta\ut\in V\tens V$. Then
\eqn{lgamma}{ g=\theta_L(\eta\uo)\theta(\eta\ut)\in\Omega^1M\tens_M
\Omega^1M\subset M\tens M\tens M.}
In practice $\eta$ will still tend to be in some completed tensor
product when we work with the universal calculus, and it should be
non-degenerate in some sense. As soon as we pass to nonuniversal
calculi, $V$ will tend to be finite-dimensional and these
subtleties will not arise. Therefore, for simplicity, we stress the
version with $\gamma$. Finally, we have seen that vanishing of the
torsion and cotorsion classically means metric compatibility in a
skew symmetrized or `differential form' sense. However,  quantum
differential forms in the universal calculus are not skew
symmetrized  so at least in this setting we can reasonably take
$D\gamma=0$ and $D\theta=0$ as the correct conditions for `torsion
free and metric compatible'.

\begin{propos} We define the cotorsion form $\Gamma\in\Omega^2M
\tens_M\Omega^1M$ as
corresponding to $D\gamma$ under $\id\tens s_\theta$.  It coincides
with $(T_\gamma\tens\id)g$ where $T_\gamma$ is the right $M$-module
map corresponding in the dual frame resolution to the torsion of
$\gamma$. Moreover,
\[ \Gamma=\extd g+ (\id\tens T)g.\]
\end{propos}
\proof  The explicit computations are similar (with a left-right
reversal) to those
already made for $\bar D\theta$ etc., so we omit them. Comparing
the result of $\Gamma=(T_\gamma\tens\id)g$ with that of $(\id\tens
T)g$, we find that they differ by $\extd g$.
\endproof

Here $(\id\tens T)g$ is the torsion form (i.e the torsion tensor
viewed in $\Omega^1M\tens_M\Omega^2M$ via the quantum metric). So
the proposition says that the cotorsion minus the torsion is $\extd
g$. Since the universal differential calculus has trivial
cohomology, we see that for a torsion free connection, a metric
$\gamma$ has zero cotorsion {\em iff} the corresponding $g$ is closed.

As a third formulation of the quantum metric, one should be able to
identify $\Omega^1M\tens_M\Omega^1M$ with $(P\tens (V\tens V))^\qg
$, i.e. with tensorial $V\tens V$ 0-forms. This is a different
quantum generalisation of $\nabla$ from the derivation property, as
we have seen in detail in the proof of Proposition~2.7 in
Section~2. This third approach, however, will be considered
elsewhere because it appears to need nonuniversal calculi for its
proper formulation. The main results above clearly do extend in
principle to nonuniversal calculi $(M,\Omega^1(M))$. Thus, a frame
resolution means a choice of quantum principal bundle with
nonuniversal calculus $(P,\Omega^1(P),\qg ,\Omega^1(\qg
),V,\theta)$ as in \cite{BrzMa:gau}\cite{BrzMa:dif} and a right
strongly tensorial form $\theta:V\to\Omega^1(P)$ such that the
induced map $s_\theta:\CE \to\Omega^1(M)$ is an isomorphism.
Similarly, we may introduce a metric as   $\gamma:V^*\to
\Omega^1(P)$ such that the induced $s_\gamma$ is an isomorphism. The
underlying theory of associated bundles with
nonuniversal calculi needs to be developed first in order to
proceed further.

Finally, we recall that a quantum principal bundle $(P,\qg)$ is
called trivial\cite{BrzMa:gau} if there is a convolution-invertible
unital linear map $\Phi:\qg \to P$ which intertwines the right
regular coaction of $\qg $ on itself and $\Delta_R$ on $P$. In
terms of extension theory of algebras, one says that the extension
is cleft, and one knows that $P$ is then a cocycle cross product
algebra. A strong connection $\omega$ in this case is equivalent to
a `gauge field' $A:\qg\to\Omega^1M$ such that $A(1)=0$. We
have\cite{BrzMa:gau}
\eqn{trivcon}{ \omega=\Phi^{-1}*A*\Phi+\Phi^{-1}*\extd\Phi}
where $*$ is the convolution product $x*y=\cdot\circ(x\tens
y)\circ\Delta$ for maps $x,y$ from $\qg$. From \cite{BrzMa:gau}, we
also know that left strongly tensorial forms such as (in our
present case) $\gamma:V^*\to (\Omega^1M)P$ on a trivial bundle are
in 1-1 correspondence with `matter fields' $\cvb:V^*\to\Omega^1M$
via $\gamma=\cvb*_R\Phi$. (Similarly, $\theta_L=\vb_L*_R\Phi$ for
some $\vb:V\to\Omega^1M$ in the alternative formulation of the
quantum metric.) Here we extend the $*$ notation to the convolution
right action with respect to $\Delta_R$ in place of $\Delta$. In a
similar way, one finds that right strongly tensorial forms such as
(in our case) $\theta:V\to P\Omega^1M$ are in 1-1 correspondence
with `matter fields' $\vb:V\to\Omega^1M$ via
$\theta=\Phi^{-1}*_L\vb$, where we now also use $*_L$ for the
convolution left action with respect to $\Delta_L$. Moreover, that
the maps $s_\theta$ and $s_\gamma$ are invertible correspond
respectively to invertibility of the maps
\eqn{svb}{s_\vb:M\tens
V\to\Omega^1M, \quad s_\vb(m\tens v)=m\vb(v),\quad s_\cvb: V^*\tens
M\to\Omega^1M,\quad s_\cvb(w\tens m)=\cvb(w)m} for $m\in M,v\in
V,w\in V^*$. The first case follows easily from the description in
\cite{BrzMa:gau} of the associated bundle $\CE $ in the trivial
bundle case as $M\tens V\isom \CE $, combined with Lemma~3.1, while
the second is similar. We similarly need $s_{\vb_L}$ an isomorphism
in that setting. We call these particular `matter fields' $\vb$ and
$\cvb$ the {\em quantum $V$-bein or $V$-cobein} respectively, and
$\vb_L$ a left-handed quantum $V$-bein. They are global
parallelisations of $\Omega^1M$ as a basis of 1-forms over $M$
acting either from the left or from the right.

\begin{propos} For a trivial quantum principal bundle frame resolution
in terms of gauge fields
and $V$-(co)beins, the covariant derivative is
\[ \nabla=1\tens \id-(\id\tens \vb)s_\vb^{-1}
-s^{-1}_{\vb\alpha}\cdot (A*_L\vb)(s^{-1\alpha}_\vb),\]
where $s_{\vb\alpha}^{-1}\tens s_\vb^{-1\alpha}$ denotes the output
of $s^{-1}_\vb$ acting on $\Omega^1M$. The torsion and cotorsion
correspond to
\[\bar D_A\vb=\extd\vb+A*_L\vb,\quad D_A\cvb=\extd\cvb+\cvb*_RA\]
and the metric has the form $g=\cvb(f^a)\vb(e_a)$, where $\{e_a\}$
is a basis of $V$ and $\{f^a\}$ is a dual basis.
\end{propos}
\proof  This follows from Proposition~2.3 and the form of $s_\vb^{-1}$
deduced from
(\ref{svb}), and the isomorphism $M\tens V\isom \CE $ from
\cite{BrzMa:gau}. The torsion tensors and cotorsion tensor
correspond to $\bar D\theta=\Phi^{-1}*_L D_A\vb$,
$D\gamma=(D_A\cvb)*_R\Phi$, for some $\bar
D_A\vb:V\to\Omega^2M,D_A\cvb:V^*\to\Omega^2M$, since they are again
right and left strongly tensorial. One immediately finds them as
shown. One similarly has $D_A\vb_L=\extd\vb_L+\vb_L*_RA$ in that
setting for the quantum metric. \endproof

Note that the requirements for $s_\vb,s_\cvb$ (or $s_{\vb_L}$) to
be isomorphisms is the same as saying that $(M,k,V,\vb)$ and
$(M,k,V^*,\cvb)$ (or $(M,k,V,\vb_L)$) are  frame resolutions   with
trivial quantum group $\qg=k$ and $P=M$. So $M$ has a trivial
principal bundle frame resolution and/or dual frame resolution {\em
iff} it has a trivial one with $\qg=k$. On the other hand,
extending $P$ to some larger trivial frame bundle with nontrivial
structure group $\qg$ allows for a larger range of covariant
derivatives induced by different gauge fields $A$. Moreover, this
is also the `local picture' when a nontrivial bundle is glued by
patching together trivial bundles as in \cite{BrzMa:gau}.

\section{Basic constructions: $q$-homogeneous spaces and bosonisation}

In this section we show how the formalism above includes various
well-known quantum spaces arising in the theory of quantum groups
and braided groups. We provide some basic general classes of
examples as well as concrete cases such as the quantum sphere
$S_q^2$ and the quantum planes $\R^n_q$. We mainly establish the
existence of the frame resolution and  the general form of the
covariant derivative $\nabla$  for our examples, and in some cases
we obtain a quantum metric $g$.

We start with the simplest of all examples, namely $M=H$ a Hopf
algebra. As we might expect, its universal differential calculus is
`parallelisable' in the sense that it can be resolved with trivial
quantum group in the frame resolution.

\begin{propos} Let $M=H$, a Hopf algebra. Then ($H,\Omega^1H$)
has a quantum frame resolution where $P=H$ and $\qg =k$, the ground
field, and
\[ V=\ker\eps\subset H,\quad \theta(v)=Sv\o\tens v\t.\]
The covariant derivative and torsion are
\[ \nabla(h\tens g)=1\tens h\tens g-hg\o\tens Sg\t\tens g\th,
\quad T(h\tens g)=h\tens g\tens 1 -h\tens1\tens g+hg\o\tens Sg\t\tens g\th\]
extended linearly and restricted to $\Omega^1H$.
\end{propos}
\proof Here $\CE =H\tens \ker\eps$ and $s_\theta(h\tens v)=hSv\o\tens v\t$ is
an isomorphism $s_\theta:\CE \to \Omega^1H$. This map is in fact
the inverse of the well-known isomorphism $\Omega^1H\to
H\tens\ker\eps$ given by $h\tens g\to hg\o\tens g\t$. Hence this
choice of $P,V,\theta$ indeed provides a frame resolution of
$H,\Omega^1H$. On the other hand, $\ker\eps\subset \qg $ is zero so
only $\omega=0$ is possible. Then $\nabla$, $T$ necessarily have
the form stated.
\endproof

This provides a `quantum geometrical' picture of the isomorphism
$\Omega^1H\isom H\tens\ker\eps$ playing a fundamental role in the
theory of differential calculi on quantum groups. Similarly, any
left-covariant $\Omega^1(H)$ has the form $H\tens
\ker\eps/Q\isom \Omega^1(H)$ where $Q$ is a right ideal in
$\ker\eps$\cite{Wor:dif}. One can view this, as above, as coming
from a frame resolution where $V=\ker\eps/Q$ and $\theta$ is
inherited from the formula above.

The trivial frame resolution here induces only one covariant
derivative. On the other hand, we can view $\theta$ as a quantum
$V$-bein as part of any trivial quantum group principal bundle with
quantum group $\qg'$ coacting on $V$, giving a larger range of
induced covariant derivatives according to gauge fields $A$ (see
Proposition~3.7). Moreover, by an evident left-right symmetry, we
also have a left frame resolution
\eqn{thetaL}{\theta_L(v)=v\o\tens Sv\t} and
hence a metric if are given a nondegenerate $\eta\in V\tens V$. In
the finite-dimensional case when $\eta$ is viewed as  map $V^*\to
V$ by $\eta(w)=\eta(\ ,w)$, we also have a dual frame resolution
$\gamma(w)=\eta(w)\o\tens S\eta(w)\t$. The corresponding quantum
metric in either case is
\eqn{metricK}{g=\eta\uo\o\tens S(\eta\ut\o \eta\uo\t)\tens \eta\ut\t}
where $\eta=\eta\uo\tens \eta\ut$. If $\eta$ is $\qg'$-invariant
then we can view $\gamma$ as a quantum $V$-cobein on the trivial
quantum principal bundle with structure quantum group $\qg'$.

\begin{propos} Let $\qg$ be a Hopf algebra. Then $(\qg,\Omega^1\qg)$
has a frame resolution
by $P=\qg\tens \qg$, quantum $V$-bein and covariant derivative
\[ V=\ker\eps,\quad \Delta_L=\Delta-\id\tens 1,\quad \vb(v)=Sv\o\tens v\t\]
\[ \nabla(h\tens g)=1\tens h\tens g-hg\o\tens Sg\t\tens g\th-hg\o
A(g\t)S g\th\tens g\fo+hg\o A(g\t)\tens 1.\] There is a unique
gauge field $A(h)=Sh\o\tens h\t-\eps(h)1\tens 1$ with zero torsion.
For any $\Lambda\in H$,
\[ \eta=S\Lambda\o\tens\Lambda\t-S\Lambda\tens 1-1\tens\Lambda
+\eps(\Lambda)1\tens 1\in V\tens V\]
is $\qg$-invariant, and when non-degenerate it defines a quantum
metric
\[ g=S\Lambda\t\tens S^2\Lambda\o S\Lambda\th\tens \Lambda\fo
-S\Lambda\t\tens S^2\Lambda\o\tens 1
-1\tens S\Lambda\o\tens\Lambda\t+\eps(\Lambda)1\tens1\tens 1.\]
There is a unique gauge field $A(h)=h\o\tens Sh\t-\eps(h)1\tens 1$
with zero cotorsion. Both gauge fields have zero curvature.
\end{propos}
\proof We use Proposition~3.7 as explained above, regarding
$P=H\tens H$ as a trivial
bundle with structure group $H$ and $\Delta_R(h\tens g)=h\tens
g\o\tens g\t$, and regarding $\theta,\theta_L$ in Proposition~4.1
as $V$-bein and left $V$-bein. We equip $V$ with the coaction
$\Delta_L$ as stated, and easily verify that $\eta\in V\tens V$ is
invariant under the tensor product coaction. Also, both gauge
fields obey $\extd A+A*A=0$, which we can interpret as $\bar
D_A\vb=0$ in the first case and $D_A\vb_L=0$ in the second. By
applying $s_\vb$ and $s_{\vb_L}$ to these equations one knows that
$A$ composed with $\Delta_L$ or $\Delta_R$ on $V$ is fully
determined. By applying $\eps$, this determines $A$ in the two
cases uniquely. In the finite-dimensional case it makes sense to
require non-degeneracy as an isomorphism $\eta:V^*\to V$ or
equivalently as a map $V^*\tens V^*\to k$. Equivalently, we
identify $V^*=\ker\eps\subset H^*$ and require that
$\<\Lambda,(Sw)x\>$ is non-degenerate as a bilinear form on $w,x\in
V^*$.
\endproof

The nondegeneracy condition here holds, for example, when $H$ is
finite-dimensional and $\Lambda$ is a normalised integral in $H$
and $\Lambda^*$ a normalised integral in $H^*$ such that
$\<\Lambda^*,\Lambda\>$ is invertible. For then
$\eta(w)=S\Lambda\o\<\Lambda\t,w\>-1\<\Lambda,w\>$ is the 
Fourier transform on $H^*$
(restricted to $\ker\eps$ and projected to $\ker\eps$) and is
invertible, cf\cite[Cor.~1.5.6]{Ma:book}. For example, we may
certainly take the functions $\qg=\C(G)$ on a finite group $G$ and
$\Lambda$ the Kronecker $\delta$-function at the identity.

\subsection{Quantum principal homogeneous spaces}

A quantum homogeneous principal bundle is \cite{BrzMa:gau} a Hopf
algebra surjection $\pi:P\to \qg $ such that
$\Delta_R=(\id\tens\pi)\circ\Delta$ makes $P$ a quantum principal
bundle over $M=P^\qg $. A sufficient condition is that the product
map $\ker\eps|_M\tens P\to\ker\pi$ is a surjection\cite{BrzMa:gau}.
$M$ is called a principal quantum homogeneous space
cf\cite{Sch:pri}. A linear splitting $i:\qg \to P$ of $\pi$ such
that $(i\tens\id)\Ad=(\id\tens\pi)\Ad\circ i$ defines a connection
$\omega(h)=(Si(h)\o)\extd i(h)\t$, see\cite{BrzMa:gau}, called the
{\em canonical connection} associated to a splitting.  Necessary
and sufficient conditions for the canonical connection to be strong
are in \cite{HajMa:pro} and amount to $i$ a unital bicovariant
splitting.

\begin{propos}  Let $M$ be a quantum principal homogeneous space
associated to $\pi:P\to \qg $.
Then $(M,\Omega^1M)$ has a quantum frame resolution $(P,\qg
,V,\theta)$, where
\[ V=\ker\eps\cap M,\quad \theta(v)=Sv\o\tens v\t,\quad \Delta_R(v)
=v\t\tens\pi(Sv\o)\]
Let $i$ be a bicovariant unital splitting. Then the associated
canonical covariant derivative  is the restriction to $\Omega^1M$
of
\[\nabla(m\tens n)=1\tens m\tens n-m\o n\o Si\circ\pi(m\t n\t)\o
\tens i\circ\pi(m\t n\t)\t Sn\th\tens n\fo.\]
\end{propos}
\proof We first check that $\theta(v)\in P\tens M$ by computing
$(\id\tens\Delta_R)\theta(v)=Sv\o \tens v\t\tens \pi(v\th)
=Sv\o\o\tens v\o\t\tens\pi(v\t)=Sv\o\tens v\t\tens 1=\theta(v)
\tens 1$, using
coassociativity and that $V\subset M$.  We also verify that
$\Delta_R$ as stated makes $V$ a right comodule. Note that this is
such that the corresponding left handed coaction as in Section~3 is
$\Delta_L(v)=\pi(\o)\tens v\t$. This makes it clear that $\Delta_R$
is indeed a coaction (since $\Delta_L$ clearly is) and that
$\Delta_R(v)\in V\tens \qg $. Indeed, clearly $\Delta_L(v)\in
\qg \tens M$ by the same argument as for $\theta_V$, and
$(\id\tens\eps)\Delta_L(v)=\pi(v)=\eps(v\o)\pi(v\t)=\eps(v)\tens
1=0$ since $v\in M$ and then $v\in\ker\eps$. We similarly verify
that $\theta$ is an intertwiner. Putting the output of $\Delta_R$
to the far right, we have
$(\Delta_R\tens\id)\theta(v)=(Sv\o)\o\tens v\t\tens
\pi((Sv\o)\t)=Sv\t\tens v\th\tens S\pi(v\o)=Sv\bo\o\tens
v\bo\t\tens v\bt=\theta(v\bo)\tens v\bt$ by coassociativity and the
form of $\Delta_Rv=v\bo\tens v\bt$. Hence the various maps are
defined as stated. By Lemma~3.1 we have an induced $s_\theta(p\tens
v)=pSv\o\tens v\t$ and we verify that it is an isomorphism
$s_\theta:P\tens V\isom\Omega^1M$. Indeed, we define the inverse as
the restriction to $\Omega^1M$ of $s_\theta^{-1}(m\tens
n)=mn\o\tens n\t$. This has its right hand output in $M$ by the
same coassociativity argument as above. Moreover,
$(\id\tens\eps)s_\theta^{-1}(m\tens n)=mn$ so $\Omega^1M$ maps to
$P\tens V$ as required. That the two maps $s_\theta$ and
$s_\theta^{-1}$ are mutually inverse is the same elementary
computation as in Proposition~4.1. Indeed, these maps are
restrictions of the corresponding maps for $P$ as a Hopf algebra
with its trivial frame resolution. Finally, putting in the form of
$\omega$ into Proposition~3.3 immediately gives $\nabla$ as shown.
We note that this simplifies slightly on exact forms, as
\eqn{nabexa}{\nabla(\extd m)=1\tens 1\tens m-1\tens m\tens 1+m\tens
1\tens 1-m\o Si\circ\pi(m\t)\o
\tens i\circ\pi(m\t)\t Sm\th\tens m\fo}
for all $m\in M$. Also, $T(\extd m)=-\nabla(\extd m)$. Since
$\nabla$ is left derivation and the torsion tensor is a left-module
map, they are fully defined by their values on exact forms. \endproof

The most well-known nontrivial example of a principal quantum
homogeneous space is the quantum sphere $M=S_q^2$, where
$P=SO_q(3)$ as the even subalgebra of $SU_q(2)$ with usual
generators $\alpha,\beta,\gamma,\delta$, and $\qg =k[z,z^{-1}]$
with projection and induced $\Delta_R$
\[ \pi\pmatrix{\alpha&\beta\cr \gamma&\delta}
=\pmatrix{z^\h&0\cr0&z^{-\h}},\quad
\Delta_R\pmatrix{\alpha&\beta\cr \gamma&\delta}=\pmatrix{\alpha\tens
z^\h&\beta\tens z^{-\h}
\cr \gamma\tens z^\h&\delta\tens z^{-\h}},\]
restricted to $SO_q(3)$. Similarly to \cite{BrzMa:gau}, one can take
 $i(z^n)=\alpha^{2n}$ and $i(z^{-n})=\delta^{2n}$  and
verify that one has a strong canonical connection (the charge 2
monopole). One may also take a slightly more complicated $i$ in
trivial bundle `patches' if one wants closer contact with the
classical formulae\cite{BrzMa:gau}. The $S_q^2$ is the subalgebra
generated by $1, b_-=\alpha\beta, b_+=\gamma\delta,
b_3=\alpha\delta$ and is an example of the family in \cite{Pod:sph}. 
One may then compute the covariant derivatives
$\nabla (\extd b_3),
\nabla (\extd b_\pm)$ and verify that they are non zero.

\subsection{Quantum planes and other braided groups}

Other natural `quantum geometries' (in the sense of being
associated naturally with quantum group symmetries) are braided
groups $B$. These (for our purposes here) are covariant objects
under some background strict quantum group $\qg$ which, like the
$\Z_2$ of supersymmetry, induces `braid statistics' on $B$. Basic
examples (all due to the author) are quantum planes, $q$-Minkowski
space and versions $BG_q$ for all the standard quantum groups, see
\cite{Ma:introp}\cite{Ma:book}. In this section we consider quantum
group Riemannian geometry on such objects, i.e. we take $M=B$.

First of all, the braided version of Proposition~4.1 follows in
just the same way: $V=\ker\und\eps\subset B$ and $\theta(v)=\und S
v\Bo\tens v\Bt\in \Omega^1B$, where $\und\eps,\und S,\und\Delta$
are the braided group counit, antipode and coproduct (the
underlines are to remind us that braided groups are not quantum
groups in the usual way, having braid statistics). Thus every
braided group $B$ is `parallelisable' with frame resolution by
trivial quantum (or braided) group $k$. We similarly have
$\gamma=\theta_L\circ \eta$ given any isomorphism $\eta:V^*\to V$,
and $\theta_L(v)=v\Bo\tens \und S v\Bt$ as a left-handed version of
$\theta$. The corresponding quantum metric is shown in Figure~6(d)
in the appendix.

On the other hand, we know from  bosonisation
theory\cite[Thm~9.4.12]{Ma:book} that every braided group has an
equivalent quantum group $B\lbiprod
\qg$ given by adjoining the background covariance quantum group. So
we can view the trivial resolution instead
as a quantum $V$-bein.

\begin{propos} Let $M=B$ be a braided group covariant under a dual
quasitriangular
Hopf algebra $(\qg,\CR)$ by left coaction
$\Delta_L(b)=b\tildebo\tens b\tildebt$. Then $(B,\Omega^1B)$ has a
frame resolution with trivial bundle $P=B\lbiprod\qg$,
$V=\ker\und\eps\subset B$ and quantum $V$-bein $\vb(v)=\und S
v\Bo\tens v\Bt$. The covariant derivative induced by a gauge field
$A$ is
\[ \nabla(\extd b)=1\tens 1\tens b-1\tens b\tens 1+b\tens 1\tens 1
+b\Bo A(b\Bt\tildebo b\Bth\tildebo)\und S b\Bt\tildebt\tens
b\Bth\tildebt.\] Moreover, if $\eta\in V\tens V$ is nondegenerate
and $\qg$-covariant then \[ g=\eta\uo\Bo\tens (\und S
\eta\uo\Bt)\und S \eta\ut\Bo\tens \eta\ut\t\] is a quantum metric.
\end{propos}
\proof Cf\cite{BrzMa:gau} we view $B\lbiprod H$ as a quantum principal
bundle with trivialisation $\Phi(h)=1\tens h$ and
$\Phi^{-1}=S\circ\Phi$. The quantum $V$-bein $\vb$ induces a
canonical form
\[ \theta(v)=(\Phi^{-1}*_L\vb)(v)=(1\tens Sv\tildebo)\cdot(\und S
v\tildebt\Bo\tens 1)\tens v\tildebt\Bt\tens 1)\]
where the product $\cdot$ is in $B\lbiprod H$. In our case $(1\tens
h)\cdot(b\tens 1)=h\o\la b\tens h\t=b\tildebt \tens h\t
\CR(b\tildebo\tens h\o)$ depends on the quasitriangular structure if one wants
to compute $\theta$ explicitly. Similarly, we have a left-handed
quantum $V$-bein $\vb_L(v)=v\Bo\tens \und S v\Bt$ which combined
with $\eta$ gives the metric as shown. One may view $\eta$ as a map
$V^*\to V$ by $\eta(w)=(\id\tens w)\eta$ and consider that
$\gamma=\theta_L\circ \eta$ is the quantum cobein, where
\[ \theta_L(v)=(\vb_L*_R\Phi)(v)=(v\tildebt\Bo\tens 1)\tens (\und S
v\tildebt\Bt\tens Sv\tildebo).\]
Here $\Delta_R(v)=v\tildebt\tens Sv\tildebo$ is the right action on
$V$ corresponding to $\Delta_L$ and invariance of $\eta$ ensures
that $\eta$ as a map is covariant. \endproof

As the simplest example, we consider $M=B=k[x]$ the `braided
line'\cite{Ma:csta}\cite{Ma:poi}, with background quantum group
$\qg=k[\dila,\dila^{-1}]$ and $\CR(\dila^m\tens \dila^n)=q^{mn}.$
The covariance under the coaction is $\Delta_L(x^m)=\dila^m\tens
x^m$ and corresponds to the $\Z$-grading of $k[x]$ by degree. The
bosonisation $P=k[x]\lbiprod k[\dila,\dila^{-1}]$ is the quantum
plane generated by $x,\dila$ with $\dila^{-1}$ adjoined (or
`quantum cylinder'). A gauge field is any $A:k[\dila,\dila^{-1}]\to
\Omega^1k[x]$ such that $A(1)=0$.
This means a collection of 1-forms $A(\dila^m)\in\Omega^1k[x]$ for
$m\ne0$. The computations are easily made from the preceding
proposition. For example,
\[ \vb(x^m)=\sum_{r=0}^m\left[{m\atop r}\right]_q (-1)^rx^{r(r-1)\over 2}
\tens x^{m-r}\]
where $\left[{m\atop r}\right]_q$ are the q-binomial coefficients.
One may similarly compute the covariant derivative using the
$q$-trinomial coefficients for the coefficients of
$(\id\tens\und\Delta)\circ\und\Delta x^m$. For the lowest
generators, one has
\[ \nabla \extd x=A(\dila)\extd x,\quad \nabla\extd x^2
=A(\dila^2)\extd x^2+(1+q)\left((\extd x)\extd x+x A(\dila)\extd x
- A(\dila^2)x\extd x\right).\]

The $sl_n$ quantum-braided planes $M=\R^n_q$ are similarly braided
groups\cite{Ma:poi} covariant under $\widetilde{SL_q(n)}=GL_q(n)$.
Their bosonisations $\R^n_q\lbiprod GL_q(n)$ are therefore now to
be regarded via Proposition~4.4 as the {\em linear frame bundles}
of the quantum planes. In general, and unlike the classical
situation, one has different q-deformed versions of   $\R_q^n$ of
various covariance types e.g. associated to covariance under all
the dilatonic extension $\widetilde{G_q}$ of the standard matrix
quantum groups $G_q$. The general construction is
provided\cite{Ma:poi} by the theory of linear braided groups
$B=V_L(R',R)$ where $R',R$ are certain `R-matrix' data and
covariance  is under the dilatonic extension of a quantum group
obtained from the quantum matrix bialgebras\cite{frt:lie} $A(R)$. The bosonisation of these quantum braided planes
provides the construction of inhomogeneous quantum groups
$\R^n_q\lbiprod
\widetilde{G_q}$ etc. in \cite{Ma:poi}, which we understand now as
frame resolutions
by $\qg=\widetilde{G_q}$ of these various quantum planes.

The obvious case, using the $so_n$ series R-matrix, is the quantum
Poincar\'e algebra $\R_q^n\lbiprod\widetilde{SO}_q(n)$ which we
understand now as the dilaton-extended orthogonal frame bundle of
$\R_q^n$. We need the more general theory of frame resolutions,
however, to accommodate the other versions of $\R_q^n$ associated
to other quantum groups. One similarly has Minkowski versions
$\R_q^{1,3}\lbiprod \widetilde{SO}_q(1,3)$ in \cite{Ma:poi}. This
also has a spinorial version where $\R_q^{1,3}=M_q(2)$ the space of
$2\times 2$ braided hermitian
matrices\cite{Ma:exa}\cite{Ma:mec}\cite{Mey:new}; their spinorial
bosonisation $P=M_q(2)\lbiprod
\widetilde{SU_q(2)\dcross SU_q(2)}$ is computed explicitly in
\cite{MaMey:bra}
and can be viewed as a double cover of the $SO_q(1,3)$ frame
resolution. Finally, these Euclidean and Minkowski space braided
groups have known quantum metrics. A quantum metric in this context
of linear braided groups is defined\cite[Def. 10.2.14]{Ma:book} as
an isomorphism the mutually dual linear
 braided groups $V_L^*(R',R)\isom V_L(R',R)$ induced by a linear
 isomorphism $\eta$
 of the mutually dual generating vector spaces. Here (in our present
 conventions)
 the evaluation
map $\ev:V_L(R',R)\tens V_L^*(R',R)\to k$  is provided by the
braided $R$-differentiation operators\cite{Ma:fre}. The
coevaluation for this is the appropriate braided-exponential
$\exp_R(p|x)$ as a powerseries in $V_L^*(R',R)\tens V_L(R',R)$. We
refer to \cite[Chapter~10]{Ma:book} for an introduction to this
`braided analysis'. Projecting to $V=\ker\eps$ (polynomials in the
generators with no constant terms), we have an induced bilinear
form
\[ \tilde\eta=\exp_R(\eta(p)|x)-1\]
as a formal powerseries in $V\tens V$. This $\tilde\eta$ in
Proposition~4.3 then induces a formal quantum metric $g$ in the
sense of quantum group Riemannian geometry with the universal
differential calculus. Note that these linear braided groups also
have more natural nonuniversal differential calculi of the correct
classical dimension, in which case one expects the universal $g$ to
collapse down to a q-deformation of the usual flat metric
corresponding to $\eta$.

\subsection{Other bosonisations and biproducts}

We mention here some different settings for braided groups, to
which the same formulae as in Proposition~4.4 apply. Thus, one may
have $B$ covariant as a left module under a quasitriangular Hopf
algebra $\qg,\CR$ where $\CR=\CR\uo\tens\CR\ut\in\qg\tens\qg$ obeys
the axioms in \cite{Dri}. The bosonisation has a similar form
$B\lbiprod\qg$ except that this time the algebra is the cross
product by the given action of $B$ and the coalgebra is the cross
coproduct by the coaction\cite{Ma:dou} $\Delta_L(b)=\CR\ut\tens
\CR\uo\la b$. We have the same result as in Proposition~4.4 but
with this form of coaction. Thus,
\[ \nabla(\extd b)=1\tens 1\tens b-1\tens b\tens 1+b\tens1\tens 1
+b\Bo A(\CR\ut\CR'\ut)\CR\uo\la \und S b\Bt\tens \CR'\uo\la b\Bth\]
where $\CR'$ is a second copy of $\CR$.

More generally, we can think of braided groups $B\in
{}^\qg_\qg\CM$, the category  of crossed modules or quantum double
modules associated to any Hopf algebra $\qg$ with invertible
antipode. (These can also be called Drinfeld-Radford-Yetter or
DRY-modules cf\cite{Dri}\cite{Rad:str}\cite{Yet:rep}.) Here $\qg$
both acts and coacts on $B$ in a compatible way (to form
effectively an action of Drinfeld's quantum double $D(\qg)$), and
the semidirect product and coproduct or `biproduct' $B\lbiprod H$
is a Hopf algebra. Conversely, every Hopf algebra projection
$\pi:P\to
\qg$ split by a Hopf algebra map is of this form for some braided
group $B$. See cf\cite{Rad:str}\cite{Ma:skl} (the latter paper
provided the braided group formulation of this theorem of
Radford's).

We have just the same formulae as in Proposition~4.4 for $B\lbiprod
\qg$ regarded as a frame resolution of such $B$. On the other hand,
we see that general biproducts  $B\lbiprod H$ are equivalent to
special cases of the quantum homogeneous principal bundles of
Section~4.2, namely those which are trivial and where the
trivialisation $\Phi$ is a Hopf algebra map.

\begin{propos} Let $B$ be a braided group in ${}^\qg_\qg\CM$. Then
its frame resolution by $P=B\lbiprod\qg$ and $\vb(v)=\und S
v\Bo\tens v\Bt$ can be viewed as a special case of Proposition~4.3
with projection $\pi(b\tens h)=h\und\eps(b)$.
\end{propos}
\proof We apply Proposition~4.3 to this case. Here $V=\und\eps$ since $M=B$ and
$\eps|_M=\und\eps$. Since the coproduct in $B\lbiprod \qg $ has the
semidirect coproduct form $\Delta b=b\Bo b\Bt\tildebo\tens
b\Bt\tildebt$, the left coaction $\Delta_L(b)=\pi(b\o)\tens
b\t=\pi(b\Bo b\Bt\tildebo)\tens b\Bt\tildebt=b\tildebo\tens
b\tildebt$ in the proof of Proposition~4.3 coincides with the given
left coaction $\Delta_L$ used in the construction of $B\lbiprod
\qg$. Hence, in Proposition~3.7,
\align{\vb(b)\equad&&=\pi(b\o)\theta(b\t)=b\tildebo\theta(b\tildebt)
=b\tildebo Sb\tildebt\o\tens b\tildebt\t\\
&&=b\tildebo S(b\tildebt\tildebo b\tildebt\Bt\o)\tens
b\tildebt\Bt\tildebt=b\Bo\tildebo b\Bt\tildebo S(b\Bo\tildebt
b\Bt\tildebt\tildebo)\tens b\Bt\tildebt\tildebt\\ &&=(b\Bo\tildebo
b\Bt\tildebo S(b\Bo\tildebt\tildebo b\Bt\tildebt\tildebo))\und S
b\Bo\tildebt\tildebt\tens b\Bt\tildebt\tildebt\\ &&=(b\Bo\tildebo\o
b\Bt\tildebo\o S (b\Bo\tildebo\t b\Bt\tildebo\t))\und S
b\Bo\tildebt\tens b\Bt\tildebt=\und S b\Bo\tens b\Bt,} where we use
the coproduct  of $B\lbiprod \qg $, the form of $\pi$, the form of
$\theta$ from Proposition~4.3, the coproduct again, then that the
braided coproduct is covariant under the coaction. Finally, we use
the antipode $S(bh)=(S(b\tildebo h))\und S b\tildebt$ of $B\lbiprod
\qg $ and the comodule axioms to collapse to an antipode
cancellation in $\qg$. We obtain the quantum $V$-bein in
Proposition~4.4 as required.
\endproof

Actually, among trivial quantum homogeneous quantum principal
bundles in Section~4.2, there are two natural cases, namely where
$\pi$ is split by a linear map $\Phi:\qg\to P$ (intertwining the
right coaction of $\qg$ and $\Delta_R$) such that $\Phi$ is either
an algebra or a coalgebra map. In these cases $\Phi$ is
convolution-invertible with $\Phi^{-1}=\Phi\circ S$ or
$\Phi^{-1}=S\circ\Phi$ respectively. The biproduct case where
$\Phi$ is a Hopf algebra map is merely the intersection of these
two cases. More generally, one has a theory of cocycle biproducts
where $\qg$ weakly coacts on $B$ (up to a cocycle), and a theory of
dual-cocycle biproducts where $\qg$ only weakly coacts on $B$ (up
to a dual cocycle).

\subsection{Strict quantum groups as base and the quantum double}

As a very particular case of the frame resolution of braided groups
in the preceding sections, we take $M=BG_q$ the braided group
versions of the usual quantum groups $G_q$. They have been
introduced by the author in \cite{Ma:exa} and are quotients of {\em
braided matrices} $B_L(R)$ with a matrix of generators
$\vecu=\{u^i{}_j\}$, coproduct and quadratic relations
\eqn{brmat}{ \und\Delta u^i{}_j=u^i{}_k\tens u^k{}_j,\quad
R\vecu_1R_{21}\vecu_2=\vecu_2R\vecu_1R_{21}.} Here we use a version
left-covariant under the corresponding $G_q$. On the other hand,
one knows from \cite{frt:lie} that such relations are also obeyed
by certain matrix generators of $U_q(\cg)$ so one can view $BG_q$
as certain versions of the algebras $U_q(\cg)$ (the deeper reason
for this is the braided group self-duality isomorphism $BG_q\isom
BU_q(\cg)$, see \cite{Ma:skl}). Therefore if one wants to view
$U_q(\cg)$ `up side down' as some kind of coordinate ring, this is
one way to do it and $BG_q\lbiprod G_q$ is a frame resolution for
it. Actually, $BG_q\lbiprod G_q\isom G_q\dcross G_q$ (see
\cite{Ma:book}) which is essentially isomorphic to some version of
the {\em dual} of Drinfeld's double $D(U_q(g))$. But the version
with $BG_q$ explicitly expresses this dual of the quantum double in
the form of Proposition~4.4, with $M=BG_q$ as the base of a
trivial quantum principal bundle. We note also that these same
bosonisations $BG_q\lbiprod G_q$ have been considered before in
\cite{Ma:mec}, as q-deformed Mackey quantisations of a particle
moving on $BG_q$ with generalised momentum quantum group
$U_q(\cg^*)$.

The general construction behind $BG_q$ is
transmutation\cite{Ma:bg}, which associates to any dual
quasitriangular Hopf algebra $\qg,\CR$ a braided group $\und \qg$
covariant under $\qg$ by the (in our case, left) adjoint action.

\begin{propos} Let $M=\und \qg$, the braided group version of
dual-quasitriangular Hopf
algebra $\qg$. Then $\und\qg\lbiprod\qg$ is a frame resolution of
$(\und \qg,\Omega^1\und\qg)$ by Proposition~4.4, with
\[ V=\ker\eps\subset \qg,\quad \Delta_L(v)=v\o Sv\th\tens v\t,\quad
 \vb(v)=\CR((Sv\t)\tildebo\tens v\o)(Sv\t)\tildebt\tens v\th\]
Moreover, if $\qg$ is factorisable with induced linear isomorphism
$\CQ=\CR_{21}\CR:\qg^*\to \qg$  then we have a quantum cobein
\[\cvb=\vb_L\circ
\CQ^{-1}\circ S,\quad \vb_L(v)=v\o\tens \CR((Sv\th)\tildebo\tens
v\t)(Sv\th)\tildebt\] and hence a quantum metric $g$.
\end{propos}
\proof Cf\cite{Ma:bg} but now in a left-handed form,
the structure of $\und\qg$ is
\[ \Delta_Lh=h\o Sh\th\tens h\t,\quad h\und\cdot g
=\CR(g\tildebo\tens Sh\t) h\o g\tildebt ,\quad \und S h
=\CR((Sh\t)\tildebo\tens h\o)(Sh\t)\tildebt,\]
in terms of the original Hopf algebra structure of $\qg$. The braided coproduct
and counit coincide with the coproduct and counit of $\qg$. Applying
this in Proposition~4.4 gives the formula for $\vb$. Similarly for
$\vb_L$. Also, it is well-known cf\cite{ResSem:mat} that the `inverse quantum Killing
form' $\CQ=\CR_{21}\CR$ is $\Ad$-invariant, and in braided group
theory it becomes a braided group homomorphism $\und
\qg^*\to\und\qg$ by $\CQ(w)=(\id\tens w)\CQ$. If we assume that $\CQ$ is
a linear isomorphism (the so-called factorisable case\cite{ResSem:mat}) then
we have the desired left adjoint coaction-invariant `Killing form'
$\eta=\CQ^{-1}(Sf^a)\tens e_a-1\tens 1$, where $\{e_a\}$ is a basis
of $\qg$ with dual basis $\{f^a\}$. The only subtlety is the
antipode needed for the correct coevaluation element
$\coev=Sf^a\tens e_a\in
\und\qg^*\tens\und\qg$, see \cite[Prop.~9.4.11]{Ma:book}. \endproof

The latter `factorisability' assumption applies to
finite-dimensional quantum groups, but it also
holds\cite{ResSem:mat} in a formal power-series setting (after
allowing suitable square-roots and logarithms of generators) for the
standard quantum groups such as $G_q$ and hence $BG_q$ (since this
coincides as a linear space with $G_q$). 

Finally, we outline a different frame resolution of $BG_q$, this
time as a left module under a quasitriangular Hopf algebra
$\qg=U_q(\cg)$   (as in Section~4.3). Here we view, by definition,
that $U_q(\cg)\equiv G_q^*$, i.e. we regard it `up side down' as
the $q$-deformed coordinate ring of the Drinfeld-dual group with
Lie algebra $\cg^*$. Then $BG_q\lbiprod G_q^*\isom D(U_q(g))$ (not
its dual as before). The general setting here is best covered by using $\und \qg$ the
braided version of a quasitriangular Hopf algebra $\qg,\CR$ in
\cite{Ma:bra} (not dual-quasitriangular as before).  One also has $\und H\lbiprod H\isom H\codcross H$,
see \cite{Ma:book}. and in the factorisable case one has $\und H\lbiprod
H\isom D(H)$, the Drinfeld quantum double\cite{Dri}. Moreover,
$\und \qg=\qg$ as an algebra.

\begin{propos} Let $M=H$ be a quasitriangular Hopf algebra. Then
$(\qg,\Omega^1\qg)$
has a quantum frame resolution by the `quantum double' in the form
$\und\qg\lbiprod
\qg$ and the quantum $V$-bein
\[ V=\ker\eps\subset H,\quad \Delta_L(v)=\CR\ut\tens\Ad_{\CR\uo}(v)
,\quad \vb(v)
=X\uth\cu^{-1}(Sv\t)SX\ut\tens X\uo v\o\]
where $X=\CR_{12}\CR_{13}\CR_{23}\in\qg^{\tens 3}$,
$\cu=(S\CR\ut)\CR\uo$ and $\Ad$ is the left quantum adjoint action.
The element $\eta=(S\tens
\id)(\CR_{21}\CR)-1\tens1\in V\tens V$ is $\Ad$-invariant and in the
factorisable case
induces a quantum metric via the left quantum $V$-bein
\[ \vb_L(v)=v\o X\uth\tens \cu^{-1} (SX\ut)(Sv\t)X\uo\]
\end{propos}
\proof  We use the braided coproduct and braided antipode of $\und\qg$ as\cite{Ma:bra}\cite{Ma:book}
\[ \und\Delta b=b\o S\CR\ut\tens\Ad_{\CR\uo}(b\t),\quad \und S b
=\cu^{-1}(S\CR\ut)(Sb)\CR\uo.\]
 Writing out $\Ad_h(g)=h\o g
Sh\t$ explicitly and using Drinfeld's quasitriangularity axioms,
one may compute $\vb(v)=\und S v\Bo\tens v\Bt$ and
$\vb_L(v)=v\Bo\tens\und S v\Bt$ as shown in terms of the structure
of $\qg$ (among many other ways to write these objects). This is a
straightforward Hopf-algebra calculation. Meanwhile, $\Ad$-invariance of $(S\tens\id)K$ is
well-known and given explicitly in \cite[Chapter 2]{Ma:book}. One
may go on and write the quantum metric $g=\vb_L(K\uo)\vb(K\ut)$
explicitly as a product of several copies of $\CR$.
\endproof

Thus the braided version $BU_q(\cg)$, which has the same algebra as
$U_q(\cg)$, has frame resolution $BU_q(\cg)\lbiprod U_q(\cg)$. This
can be applied to the reduced quantum group enveloping algebras at
roots of unity (which are finite-dimensional), or applied in the
formal power-series setting of \cite{Dri}. In another other version
of Proposition~4.7 we may take $BG_q$ in place of $BU_q(\cg)$ since
these are essentially isomorphic in the factorisable case. Then
$BG_q$ has a frame resolution by $BG_q\lbiprod U_q(\cg)$. Likewise,
we have a version of Proposition~4.6  where we replace $BG_q$ by a
suitable  (right $\Ad$-action covariant) version of $BU_q(\cg)$. In
this form, one may take $\eta=(\id\tens S)(\CQ)-1\tens 1\in
BU_q(\cg)\tens BU_q(\cg)$.

There are many other braided groups beyond those discussed above.
For example, it is obvious from Lusztig's book\cite{Lus} that in
his approach to the structure of $U_q(\cg)$ one effectively views
$U_q(n_+)$ as a braided group with phase-factor braid statistics as
in \cite{Ma:csta}. The above results provide a step towards a
`quantum group Riemannian geometry' of such objects as well, albeit
far removed from our original physical motivation of $q$-deforming
usual geometry.

\appendix
\section{Braided group `diagrammatic' Riemannian geometry}

Here we give the formalism developed above in a different
`diagrammatic' setting of braided group gauge
theory\cite{BrzMa:coa}\cite{Ma:diag}, i.e. where the gauge group
has braid statistics (this should not be confused with Section~4
where we gave some examples where the base $M$ was a braided
group). In particular, \cite{Ma:diag} developed principal bundles,
connections, associated bundles etc at the level of braid and
tangle diagrams, which theory we `update'  now to include the
elements of Riemannian geometry above. As well as being more
general and having potentially different examples than the quantum
group case (cf. the anyonic (or $\Z_n$-graded) gauge theory in
\cite{Ma:diag}), the diagrammatic theory provides a different style
of proofs which, for trivial braid statistics, reduces to the
quantum group gauge theory. For super ($\Z_2$-graded) gauge theory
we just take bose-fermi statistics, so the formulae for
super-quantum group Riemannian geometry can be read off from these
diagrams. In general we work in a braided category\cite{JoyStr:bra}
where, for any two objects there is a braiding
$\Psi=\epsfbox{braid.eps}$ implementing their exchange. We denote
$\Psi^{-1}=\epsfbox{braidinv.eps}$. Extending this notation,
algebra products are denoted $\epsfbox{prodfrag.eps}$ and
coproducts or coactions denoted $\epsfbox{deltafrag.eps}$. Maps are
provided by `wiring' outputs into inputs, with maps flowing
generally downwards. The unit object $\und 1$ for the tensor
product is denoted by omission. This is the `diagrammatic braided group theory'
introduced in \cite{Ma:bra}. See \cite{Ma:introp}\cite{Ma:book}.

Thus we consider an algebra $M$ in a braided category, the
differential calculus $\Omega^1M$ defined diagrammatically and a
braided group principal bundle $P,B$. Here $B$ is a braided group
or Hopf algebra with braid statistics\cite{Ma:bra}. If $V$ is a
right $B$-comodule in the braided category, we have an associated
bundle $\CE =(P\tens V)^B$ as before, where fixed points are
defined categorically as equalisers and where $P\tens V$ has the
braided tensor product coaction. This theory is in \cite{Ma:diag}.
We assume that our braided category has direct sums and appropriate
flatness properties, as explained in \cite{Ma:diag}, and adopt the
corresponding abuses of notation. To this we add:

\begin{lemma} Right strongly tensorial $\theta:V\to P\Omega^nM$ are in 1-1
correspondence with left $M$-module morphisms $\CE \to \Omega^nM$.
\end{lemma}
\proof This is shown in Figure~1 for $n=1$ (the general case looks
just the same).
We (a) apply the braided coaction $P\to P\tens B$ as shown and find
that it acts trivially, hence the morphism factors through
$\Omega^1M$. Conversely, (b) constructs $\theta:V\to P\Omega^nM$
given $s:\CE \to\Omega^nM$. The proofs follow exactly the steps for
the $n=0$ case in \cite{Ma:diag}, so we omit the detailed
diagrammatic verification. The morphism $\eta:\und 1\to B$ is plays
the role of the unit `element' of the braided group (and should not
be confused with the local metric $\eta$ in the main sections of
the paper).
\endproof
\begin{figure}
\[ \epsfbox{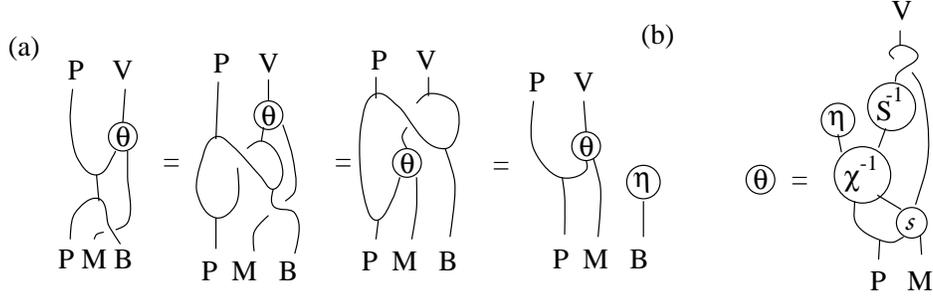}\]
\caption{Proof of Lemma~A.1}
\end{figure}

We then define a braided frame resolution as a pair $P,\theta$ such
that $P$ is a braided principal bundle and the morphism
$s_\theta=(\cdot\tens\id)(\id\tens\theta)$ is an isomorphism. We
then induce similar isomorphisms for $\Omega^1M\tens_M\Omega^1M$ by
$\id\tens s_\theta$.

\begin{propos} Given a braided frame resolution of $M$ and a
connection $\omega$,
we define $\nabla:\Omega^1M\to
\Omega^1M\tens\Omega^1M$ as the covariant derivative
$D=(\id-\Pi_\omega)\extd$ viewed under he
above isomorphisms. It is computed in Figure~2 and is a derivation
with respect to multiplication in the first factor.
\end{propos}
\proof This is shown in Figure~2, using the expression for
$\id-\Pi_\omega$ in \cite{Ma:diag}.
The derivation formula is then immediate from the final result for
$\nabla$. The unmarked $\epsfbox{deltafrag.eps}$ denotes the
coaction $\Delta_R:P\to P\tens B$.
\endproof
\begin{figure}
\[ \epsfbox{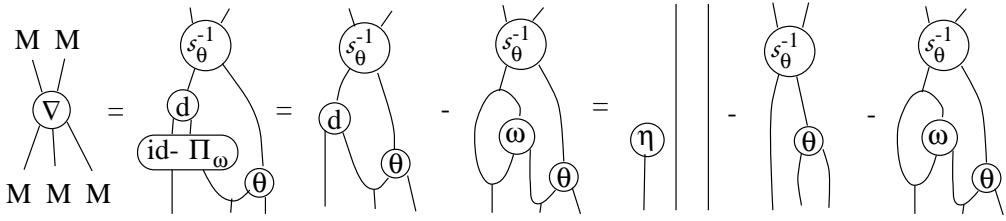}\]
\caption{Diagrammatic form of $\nabla$}
\end{figure}

Thus, we are able to proceed along the same lines as in the quantum
group case, with similar results. We let
\[ \Delta_L=(S^{-1}\tens\id)\Psi^{-1}\circ\Delta_R\]
which\cite{Ma:introp} makes $P$ into a braided left $\bar B$-module
algebra in the category with reversed braiding. Here $\bar B$
denotes $B$ with the opposite product $\cdot\circ\Psi^{-1}$ and
is\cite{Ma:introp} a braided group in the category with reversed
braiding. We define $\bar D=(\id-\bar\Pi_\omega)\extd$ on $P$,
where $\bar\Pi_\omega$ is the left-handed version (the mirror
reflection) of the diagram for $\Pi_\omega$ in \cite{Ma:diag}.

\begin{propos} The following are equivalent: (a) A connection $\omega$
is left strong. (b) The equality in Figure~3(b) holds; (c) The
equality in Figure~3(c) holds. In this case, $\bar D$ preserves
right strongly tensorial forms on $P$.
\end{propos}
\proof The condition for $\omega$ to be left strong (to preserve left
strongly tensorial forms) is given in \cite{Ma:diag}. Applying
$P\tens_M(\ )$ to this diagram and making a product in $P$, one
finds that it factors through $\chi:P\tens_MP\to P\tens B$.
Cancelling this gives the first equality in Figure~3(a) as the
condition for left strongness. The left hand side of this is,
however, equal to its mirror image (shown on the right in
Figure~3(a)). The proof of this is the lower line in Figure~3(a).
We first insert a trivial `antipode loop'. The next equality is the
comodule property of $\Delta_R$. We then use $\Ad$-invariance of
$\omega$, cancel a resulting antipode-loop and finally rearrange to
recognize the right hand side of the upper line in Figure~3(a).
Parts (b) and (c) immediately follow as equivalent to these two
versions of the left-strongness conditions: We precompose with the
coproduct of $B$, apply the antipode to the free leg thus created
and join up as a product in $B$ from the appropriate side in such a
way as to create an antipode loop cancellation in each term. As a
corollary, if $\omega$ is such that $D$ preserves left strongly
tensorial forms then $\bar D$ preserves right strongly tensorial
forms since the condition for the latter is the mirror image of the
condition for the former. We use the mirror image of the proof of
\cite[Prop. 4.2]{Ma:diag} (this part of the proof does not actually
require the Galois condition). \endproof
\begin{figure}
\[ \epsfbox{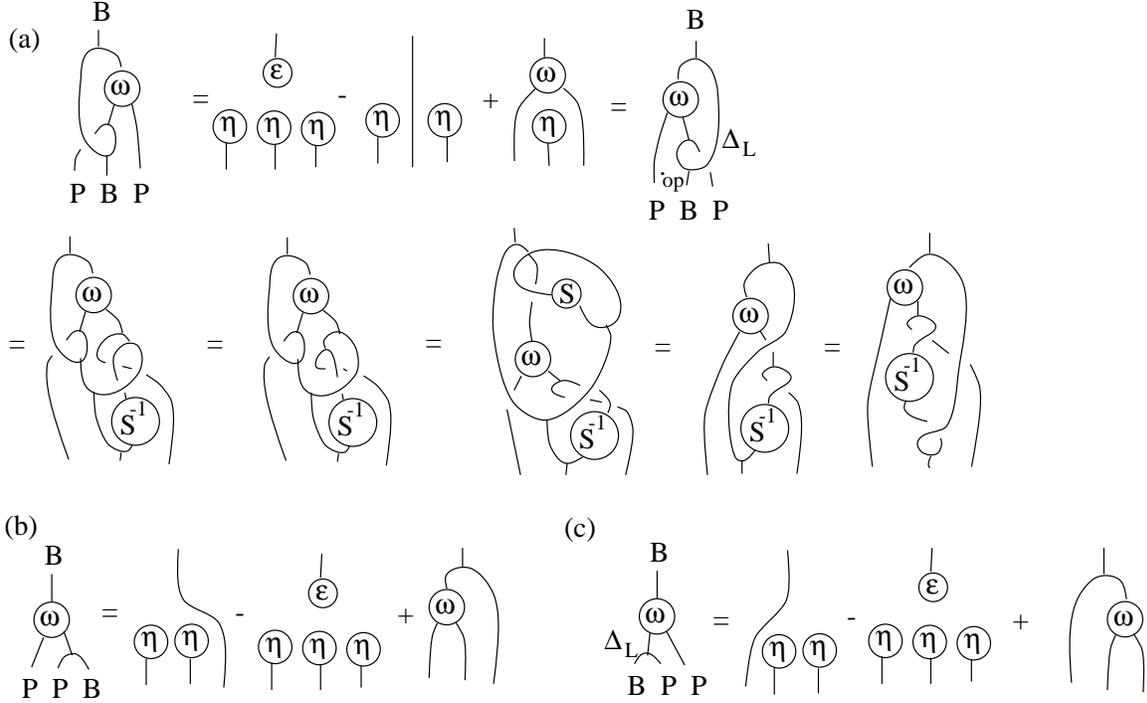}\]
\caption{Equivalent versions of the strongness condition for $\omega$}
\end{figure}

\begin{propos} We define the braided torsion tensor: $\Omega^1M\to
\Omega^1M\tens_M\Omega^1M$
as $\extd-\nabla$. This is shown in Figure~4(a). It corresponds
under the above isomorphisms to $\bar D\theta$ as a right strongly
tensorial 2-form, shown in Figure~4(b).
\end{propos}
\proof The proof follows the quantum group case. The exterior derivative
has the identical form (but written diagrammatically) and we
subtract $\nabla$. We also write the coaction on $P$ in $(P\tens
V)^B$ in $\nabla$ as a coaction on $V$ (see \cite[Prop.
4.4]{Ma:diag}). For the second part, we graft on a product with
$P\tens_M$ from the left and cancel a copy of $\chi$ through which
the diagram factors. The resulting diagram is same as grafting on a
product with $P$ from the left to $\bar D\theta$ shown in
Figure~4(b). This is what is required according to Lemma~A.1.
Finally, $\bar D\theta=(\id-\bar\Pi_\omega)\extd\theta$ when this
is computed in a similar way to Proposition~A.3.
\endproof
\begin{figure}
\[ \epsfbox{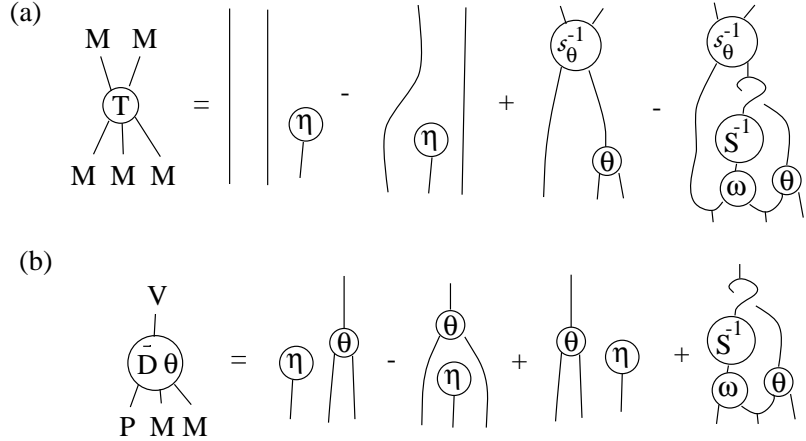}\]
\caption{Braided torsion tensor}
\end{figure}

Next, we can formulae a metric as, by definition, corresponding to
$\gamma$ a left-strongly tensorial form making $(P,B,V^*,\gamma)$ a
(left handed) frame resolution. The associated metric $g:\und
1\to\Omega^1M\tens_M\Omega^1M=\Omega^2M$ is given in Figure~5(a),
along with the associated left-handed version of the isomorphism
$s_\gamma:(V^*\tens P)^B\isom \Omega^1M$. This assumes $V$ is rigid
in the sense of a dual with coevaluation $\cap:\und 1\to V^*\tens
V$. More generally, one should work instead with $(P,B,V,\theta_L)$
a left handed frame resolution and an invariant morphism $G:\und 1\to
V\tens V$ (say). Then $g=\cdot\circ(\theta_L\tens\theta)\circ G$. The
diagrams in this case are very similar so we will not repeat them
explicitly.

\begin{propos}  The cotorsion form $\Gamma$
defined as corresponding to $D\gamma$ under $\id\tens s_\theta$
coincides with $(T_\gamma\tens\id)\circ g$, where $T_\gamma$ is the
torsion of $\gamma$ computed in the dual frame resolution evaluated
on $g$. Moreover, $\Gamma=\extd g+(\id\tens T)\circ g$.
\end{propos}
\proof The first part is shown in Figure~5(b) and follows at once from
the definition of the metric $g$ in terms of $\gamma$. Here
$T_\gamma$ is shown in the dotted box  as obtained from $D\gamma$
via $s_\gamma^{-1}$. Figure~5(c) gives the explicit form of
$D\gamma$ as a special case of \cite{Ma:diag}. Part (d) then gives
the explicit form of $\Gamma$ given by combining $D\gamma$ with
$\theta$ as in part (b). Similarly combining $\bar D\theta$ from
Figure~5 with $\gamma$ gives  an expression with similar last two
terms; comparing them we see that the difference is precisely
$\extd g$.
\endproof
\begin{figure}
\[ \epsfbox{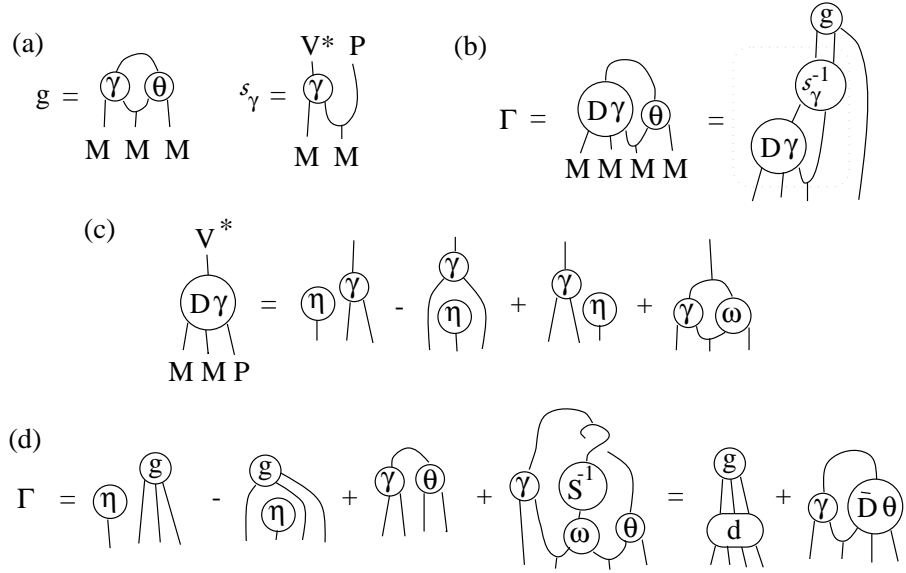}\]
\caption{Braided metric and its cotorsion}
\end{figure}

Finally, if the frame resolution bundle is trivial, with
trivialisation $\Phi:B\to P$, we define $\vb:V\to\Omega^1M$
corresponding to $\theta$ via Figure~6(a). Also shown is the
induced isomorphism $s_\vb:M\tens V\isom\Omega^1M$ by
$s_\vb=s_\theta\circ\theta_{\CE} $, where $\theta_{\CE}:M\tens
V\isom
\CE $ is given in \cite{Ma:diag}. This makes it clear that $s_\theta$ invertible
is equivalent to $s_\vb$ invertible. Similarly, the existence of a
dual frame resolution $\gamma$ is equivalent to
$\cvb:V^*\to\Omega^1M$ such that $s_\cvb$ is invertible (see
Figure~6(c)). Similarly, we already know\cite{Ma:diag} that strong
connections $\omega$ correspond to gauge fields $A:B\to\Omega^1M$
via $\omega=\Phi^{-1}*A*\Phi+\Phi^{-1}*\extd\Phi$. Putting this
into $\nabla$ from Figure~2 and moving the coaction on $P$ over to
$\Delta_L$ in $V$ (as in the preceding proof) gives the resulting
covariant derivative $\nabla$ in terms of $\vb,A$ as shown in
Figure~6(b). In fact, the resulting `local' formulae take the same
form as in the quantum group case at the end of Section~3 if we use
the convolution product notation. One has similarly that $D\gamma$
and $\bar D\theta$ correspond to $D_A\cvb$ and $\bar D_A\vb$ as in
Proposition~3.7, written diagrammatically as morphisms. The
formulae for a left-handed frame resolution $\vb_L$ are similar.
\begin{figure}
\[ \epsfbox{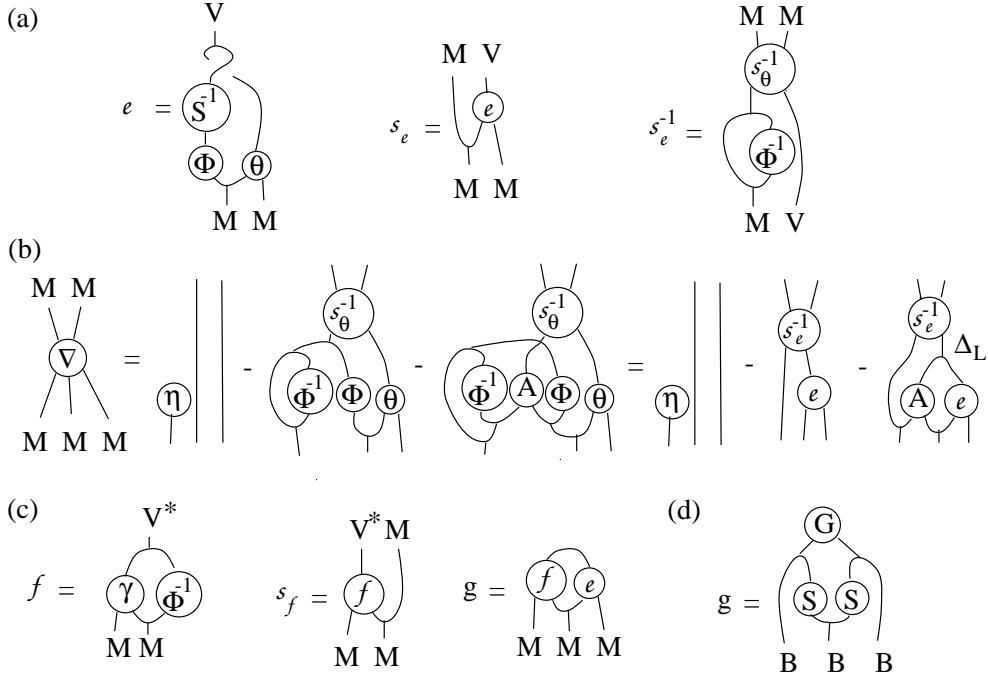}\]
\caption{Local picture in terms of braided $V$-bein and cobein,
and gauge field $A$}
\end{figure}

We see that braided theory goes through along the lines of the
quantum group case. One has braided versions of all the examples in
Section~4 as well: braided groups, braided homogeneous spaces and
braided cross products such as $\C_q^2\lcross BGL_q(2)$ in the
parallel way. For example, the braided analogue of Proposition~4.1
is $M=B$ a braided group, $V=\ker\und\eps$ (defined now as an
equaliser), and $\theta=( S\tens\id)\circ\Delta$. There is also a
left handed $\theta_L=(\id\tens S)\circ\Delta$ and hence if there
is an invariant morphism $G:\und 1\to V\tens V$ (say) we have a
braided metric as shown in Figure~6(d). The general braided theory
is, however, potentially better behaved as regards the $\Ad$ bundle
and other properties than the quantum group theory through the
imposition of natural `braided commutativity', see \cite{Ma:diag}.
Moreover, the braided setting allows one to read off the $\Z_2$ and
$\Z_n$-graded versions of the theory by inserting the relevant
braid statistics phase factor at each braid crossing.

%\bibliographystyle{unsrt}
%\bibliography{biblio}

\end{document}